\documentclass[reqno,11pt]{article}
\usepackage{psfig, amsmath, amstexnb, amsthm}
\usepackage{amssymb}  
\input stochbif.sty

\begin{document}


\title{Pathwise description of dynamic pitchfork bifurcations \\
with additive noise}
\author{Nils Berglund and Barbara Gentz}
\date{}

\maketitle

\begin{abstract}
\noindent
The slow drift (with speed $\eps$) of a parameter through a pitchfork
bifurcation point, known as the dynamic pitchfork bifurcation, is
characterized by a significant delay of the transition from the unstable to
the stable state. We describe the effect of an additive noise, of intensity
$\sigma$, by giving precise estimates on the behaviour of the individual
paths. We show that until time $\sqrt\eps$ after the bifurcation, the paths
are concentrated in a region of size $\sigma/\eps^{1/4}$ around the
bifurcating equilibrium. With high probability, they leave a neighbourhood
of this equilibrium during a time interval $[\sqrt\eps,
c\sqrt{\eps\abs{\log\sigma}}\,]$, after which they are likely to stay close
to the corresponding deterministic solution. We derive exponentially small
upper bounds for the probability of the sets of exceptional paths, with
explicit values for the exponents.
\end{abstract}

\leftline{\small{\it Date.\/} August 4, 2000.}
\leftline{\small 2000 {\it Mathematics Subject Classification.\/} 37H20,
60H10 (primary), 34E15, 93E03 (secondary).}


\section{Introduction}

Physical systems are often described by ordinary differential equations
(ODEs) of the form
\begin{equation}
\label{i1}
\dtot xs = f(x,\lambda),
\end{equation}
where $x$ is the state of the system, $\lambda$ a parameter, and $s$
denotes time. The model \eqref{i1} may however be too crude, since it
neglects all kinds of perturbations acting on the system. We are interested
here in the combined effect of two perturbations: a slow drift of the
parameter, and an additive noise.

\begin{figure}
 \centerline{\psfig{figure=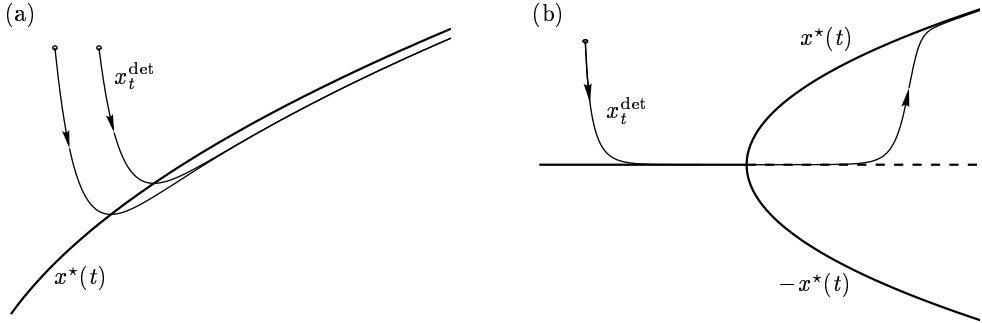,width=130mm,clip=t}}
 \caption[]
 {Solutions of the slowly time-dependent equation \eqref{i2}
 represented in the $(t,x)$-plane.  
 (a) Stable case: A stable equilibrium branch $x^\star(t)$
 attracts nearby solutions $\xdet_t$. Two solutions with different initial
 conditions are shown. They converge exponentially fast to each other, as
 well as to a neighbourhood of order $\eps$ of $x^\star(t)$. (b) Pitchfork
 bifurcation: The stable equilibrium $x=0$ becomes unstable at $t=0$
 (broken line) and expels two stable equilibrium branches $\pm
 x^\star(t)$. A solution $\xdet_t$ is shown, which is attracted by $x=0$,
 and stays close to the origin for a finite time after the bifurcation.
 This phenomenon is known as bifurcation delay.}
\label{fig1}
\end{figure}

A slowly drifting parameter $\lambda=\eps s$, (with $\eps\ll1$), may model
the deterministic change in time of some exterior influence, such as the
climate acting on an ecosystem or a magnetic field acting on a
ferromagnet.  Obviously, nontrivial dynamics can only be expected when
$\lambda$ is allowed to vary by an amount of order $1$, and thus the system
has to be considered on the time scale $\eps^{-1}$. This is usually done by
introducing the {\it slow time\/} $t=\eps s$, which transforms \eqref{i1} into the
singularly perturbed equation
\begin{equation}
\label{i2}
\eps\dtot xt = f(x,t).
\end{equation}
It is known that solutions of this system tend to stay close to stable
equilibrium branches of $f$ \cite{Grad,Tihonov}, see \figref{fig1}a. New,
and sometimes surprising phenomena occur when such an equilibrium branch
undergoes a bifurcation. These phenomena are usually called \defwd{dynamic
bifurcations} \cite{Benoit}\footnote{%
Unfortunately, the term \lq\lq dynamical bifurcation\rq\rq\ is used in a
different sense in the context of random dynamical systems, namely to
describe a  bifurcation of the family of invariant measures as opposed to a
\lq\lq phenomenological bifurcation\rq\rq, see for instance
\cite{Arnold}.}.
In the case of the Hopf bifurcation, when the equilibrium gets unstable
while expelling a stable periodic orbit, the bifurcation is substantially
delayed: solutions of \eqref{i2} track the unstable equilibrium (for a {\em
non-vanishing} time interval in the limit $\eps\to 0$) before jumping to
the limit cycle \cite{Shishkova,Neishtadt}. A similar phenomenon exists for
the dynamic pitchfork bifurcation of an equilibrium without drift, the
simplest example being $f(x,t)=tx-x^3$ (\figref{fig1}b). The delay has been
observed experimentally, for instance, in lasers \cite{ME} and in a damped
rotating pendulum \cite{BK}. 

These phenomena have the advantage of providing a genuinely dynamic point
of view for the concept of a bifurcation. Although one often says that a
bifurcation diagram (representing the asymptotic states of the system as a
function of the parameter) is obtained by varying the control parameter
$\lambda$, the impatient experimentalist taking this literally may have the
surprise to discover unstable stationary states of the system (s)he
investigates. The asymptotic state of the system \eqref{i1} with slowly
varying parameter $\lambda(\eps s)=\lambda(t)$ may depend not only on the
initial condition $(x_0,t_0)$, but also on the history of variation of the
parameter $\set{\lambda(t)}_{t\geqs t_0}$.

The perturbation of \eqref{i1} by an additive noise can be modeled by a
stochastic differential equation (SDE) of the form
\begin{equation}
\label{i3}
\6x_s = f(x_s,\lambda) \6s + \sigma \6W_s,
\end{equation}
where $W_s$ denotes the standard Wiener process,  and $\sigma$ measures the
noise intensity. A widespread approach is to analyse the probability
density of $x_s$, which satisfies the Fokker--Planck equation. In
particular, if $-f$ can be written as the gradient of a \defwd{potential
function} $F$, then there is a unique stationary density  $p(x,\lambda) =
\e^{-F(x,\lambda)/\sigma^2}/N$,  where $N$ is the normalization. This
formula shows that for small noise intensity, the stationary density is
sharply peaked around stable equilibria of $f$. 

That method has, however, two major limitations. The first one is that
the Fokker-Planck equation is difficult to solve, except in the linear and
in the  gradient case. The second limitation is more serious: the density
gives no information on correlations in time, and even when the density is
strongly localized, individual paths can perform large excursions. This is
why other approaches are important. A classical one is based on the
computation of first exit times from the neighbourhood of stable equilibria
\cite{FW,FJ}. 

The effect of bifurcations has been studied more recently by methods based
on the concept of random attractors \cite{CF1, Schmalfuss, Arnold}. In
particular, Crauel and Flandoli showed that according to their definition,
\lq\lq Additive noise destroys a pitchfork bifurcation\rq\rq\ \cite{CF2}.
The physical interpretation of random attractors is, however, not
straightforward, and alternative characterizations of stochastic
bifurcations are desirable. In the same way a slowly varying parameter
helps our understanding of bifurcations in the deterministic case, it can
provide a new point of view in the case of random dynamical systems.

Let us consider the combined effect of a slowly drifting parameter and
an additive noise on the ODE \eqref{i1}. We will focus on the case of a
pitchfork bifurcation, where  the questions {\it How does the additive
noise affect the bifurcation delay?} and {\it Where does the path go after
crossing the bifurcation point?} are of major physical interest.  The
situation of the drift term $f$ in \eqref{i3} depending explicitly on time
is considerably more difficult to solve than the autonomous case, and thus
much less understood. One can expect, however, that a slow time dependence
makes the problem accessible to perturbation theory, and that one may take
advantage of techniques developed to study singularly perturbed equations
such as \eqref{i2}. 
With $\lambda=\eps s$, Equation \eqref{i3} becomes
\begin{equation}
\label{i4}
\6x_s = f(x_s,\eps s) \6s + \sigma \6W_s.
\end{equation}
If we introduce again the slow time $t=\eps s$, the Brownian motion is
rescaled, resulting in the SDE 
\begin{equation}
\label{i5}
\6x_t = \frac1\eps f(x_t,t)\6t + \frac{\sigma}{\sqrt\eps} \6W_t.
\end{equation}
Our analysis of~\eqref{i5} is restricted to one-dimensional $x$. 
The noise intensity $\sigma$ should be considered as a function of $\eps$.
Indeed, since we now consider the equation on the time scale $\eps^{-1}$, a
constant noise intensity would lead to an infinite spreading of trajectories
as $\eps\to 0$. In the case of the pitchfork bifurcation, we will need 
to assume that $\sigma\ll\sqrt\eps$. 

Various particular cases of equation \eqref{i5} have been studied before,
from a non-rigorous point of view.  In the linear case
$f(x,\lambda)=\lambda x$, the distribution of first exit times was
investigated and compared with experiments in \cite{TM,SMC,SHA}, while  
\cite{JL} derived a formula for the last crossing of zero. In the case
$f(x,\lambda)=\lambda x-x^3$,  \cite{Gaeta} studied the dependence of the
delay on $\eps$ and $\sigma$ numerically, while \cite{Kuske} considered the
associated Fokker-Planck equation, the solution of which she approximated
by a Gaussian Ansatz. 

In the present work, we analyse \eqref{i5} for a general class of odd
functions $f(x,\lambda)$ undergoing a pitchfork bifurcation. We  use a
different approach, based on a precise control of the {\em whole paths}
$\set{x_s}_{t_0\leqs s\leqs t}$ of the process. The results thus contain
much more information than the probability density. It also turns out that
the technique we use allows to deal with nonlinearities in quite a natural
way. 
\begin{figure}
 \centerline{\psfig{figure=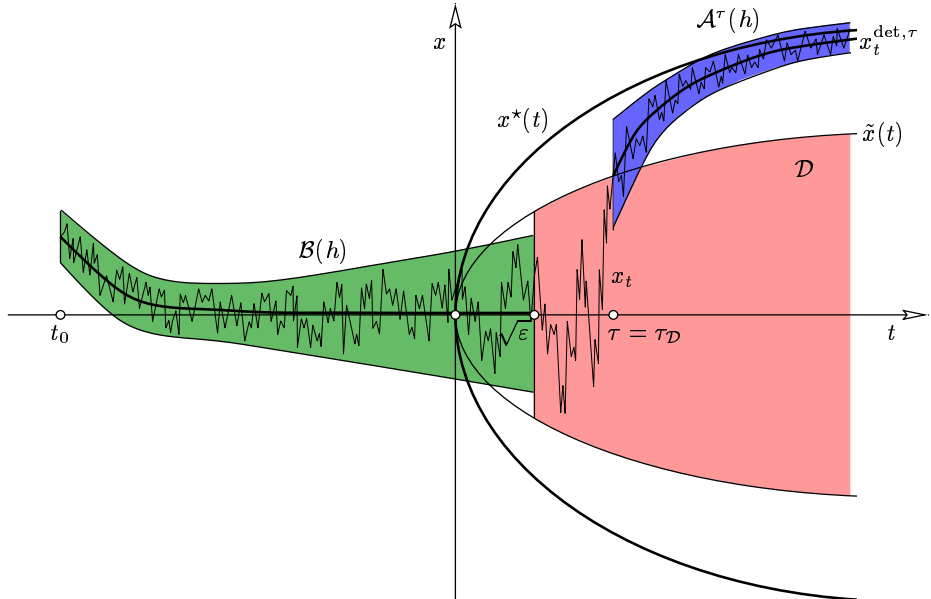,height=80mm,clip=t}}
 \caption[]
 {A typical path $x_t$ of the stochastic differential equation \eqref{i5}
 near a pitchfork bifurcation. We prove that with probability exponentially
 close to $1$, the path has the following behaviour. For $t_0\leqs
 t\leqs\sqrt\eps$, it stays in a strip $\cB(h)$ constructed around the
 deterministic solution with the same initial condition. After
 $t=\sqrt\eps$, it leaves the domain $\cD$ at a random time $\tau=\tau_\cD$,
 which is typically of the order $\sqrt{\eps\abs{\log\sigma}}$. Then it
 stays (up to times of order $1$ at least) in a strip $\cA^{\tau}(h)$
 constructed around the deterministic solution \smash{$\xdettau t$} starting at
 time $\tau$ on the boundary of $\cD$. The widths of $\cB(h)$ and
 $\cA^{\tau}(h)$ are proportional to a parameter $h$ satisfying $\sigma\ll
 h\ll\sqrt\eps$.}
\label{fig2}
\end{figure}
Our results can be summarized in the following way (see \figref{fig2}):
\begin{itemiz}
\item   Solutions of the deterministic equation \eqref{i2} starting near a
stable equilibrium branch of $f$ are known to reach a neighbourhood of
order $\eps$ of that branch in a time of order $\eps\abs{\log\eps}$. We
show that the paths of the SDE \eqref{i5} with the same initial condition are
typically concentrated in a neighbourhood of order $\sigma$ of the
deterministic solution  (Theorem~\ref{t_stable}).

\item   A particular solution of the deterministic equation \eqref{i2} is
known to exist in a neighbourhood of order $\eps$ of each unstable
equilibrium branch of $f$. Paths that start in a neighbourhood of order
$\sigma$ of this solution are likely to leave that neighbourhood in a time
of order $\eps\abs{\log\eps}$ (Theorem~\ref{t_unstable}). 

\item   When a pitchfork bifurcation occurs at $x=0$, $t=0$, the typical
paths are concentrated in a neighbourhood of order $\sigma/\eps^{1/4}$ of
the deterministic solution with the same initial condition up to time
$\sqrt\eps$ (Theorem~\ref{t_before}). 

\item   After the bifurcation point, the paths are likely to leave a
neighbourhood of order $\sqrt t$ of the unstable equilibrium before a
time $c\sqrt{\eps\abs{\log\sigma}}$ (Theorem~\ref{t_escape}). 

\item   Once they have left this neighbourhood, the paths remain with high
probability in a region of size $\sigma/\sqrt t$ around the corresponding
deterministic solution, which approaches a stable equilibrium branch of $f$
like $\eps/t^{3/2}$ (Theorem~\ref{t_approach}). 
\end{itemiz}
These results show that the bifurcation delay, which is observed in the
dynamical system~\eqref{i2}, is destroyed by additive noise as soon as the
noise is not exponentially small. Do they mean that the dynamic bifurcation
itself is destroyed by additive noise? This is mainly a matter of
definition. On one hand, we will see that independently of the initial
condition, the probability of reaching the upper, rather than the lower
branch emerging from the bifurcation point, is close to $\frac12$. The
asymptotic state is thus selected by the noise, and not by the initial
condition. Hence, the bifurcation is destroyed in the sense of \cite{CF2}.
On the other hand, individual paths are concentrated near the stable
equilibrium branches of $f$, which means that the bifurcation diagram will
be made visible by the noise, much more so than in the deterministic case.
So we do observe a qualitative change in behaviour when $\lambda$ changes
its sign,  which can be considered as a bifurcation.   

The precise statements and a discussion of their consequences are given in
Section \ref{sec_results}. In Section~\ref{ssec_rn}, we analyse the motion
near equilibrium branches away from bifurcation points. The actual
pitchfork bifurcation is discussed in Section~\ref{ssec_rb}. A few
consequences are derived in Section~\ref{ssec_rd}. Section~\ref{sec_nobif}
contains the proofs of the first two theorems on the motion near
nonbifurcating equilibria, while the proofs of the last three theorems on
the pitchfork bifurcation are given in Section~\ref{sec_pitchfork}.

\subsubsection*{Acknowledgements:}  
It's a great pleasure to thank Anton Bovier for sharing our enthusiasm. We
enjoyed lively discussions and his constant interest in the progress of our
work. The central ideas were developed during mutual visits in Berlin
resp.\ Atlanta. N.\,B. thanks the WIAS and B.\,G. thanks Turgay Uzer and
the School of Physics at Georgia Tech for their kind hospitality. N.\,B.
was partially supported by the Fonds National Suisse de la Recherche
Scientifique, and by the Nonlinear Control Network of the European
Community, Grant ERB FMRXCT--970137.


\section{Statement of results}
\label{sec_results}


\subsection{Preliminaries}
\label{ssec_rp}
We consider nonlinear It\^o SDEs of the form
\begin{equation}
\label{SDE}
\6x_t = \frac{1}{\eps} f(x_t,t) \6t + \frac{\sigma}{\sqrt{\eps}} \6W_t,
\qquad x_{t_0}=x_0,
\end{equation}
where $\{W_t\}_{t\geqs t_0}$ is the standard Wiener process on some
probability space $(\Omega, \cF, \fP)$. Initial conditions $x_0$ are
always assumed to be square-integrable with respect to $\fP$ and
independent of $\{W_t\}_{t\geqs t_0}$. All stochastic integrals are
considered as It\^o integrals, but note that It\^o and Stratonovich
integrals agree for integrands depending only on time and
$\omega$. Without further mentioning we always assume that $f$
satisfies the usual (local) Lipschitz and bounded-growth conditions
which guarantee existence and (pathwise) uniqueness of a
(strong) solution $\{x_t\}_t$ of~\eqref{SDE}. Under these conditions,
there exists a continuous version of $\{x_t\}_t$. Therefore we
may assume that the paths $\omega\mapsto x_t(\omega)$ are continuous
for $\fP$-almost all $\omega\in \Omega$. 

We introduce the notation $\fP^{\mskip1.5mu t_0,x_0}$ for the law of
the process $\{x_t\}_{t\geqs t_0}$, starting in $x_0$ at time $t_0$, and use
$\E^{\mskip1.5mu t_0,x_0}$ to denote expectations with respect to
$\fP^{\mskip1.5mu t_0,x_0}$. Note that the stochastic process $\{x_t\}_{t\geqs
t_0}$ is an (inhomogeneous) Markov process. We are interested in first
exit times of $x_t$ from space--time sets. Let $\cA \subset
\R\times[t_0,t_1]$ be Borel-measurable. Assuming that $\cA$ contains
$(x_0,t_0)$, we define the first exit time of $(x_t,t)$ from $\cA$ by 
\begin{equation}
\label{pre1}
\tau_{\cA}=\inf\bigsetsuch{t\in[t_0,t_1]}{(x_t,t)\not\in\cA},
\end{equation}
and agree to set $\tau_{\cA}(\omega)=\infty$ for those
$\omega\in\Omega$ which satisfy $(x_t(\omega),t)\in\cA$ for all
$t\in[t_0,t_1]$. For convenience, we shall call $\tau_\cA$ the {\it first
exit time of $x_t$ from $\cA$}. Typically, we will consider sets of
the form $\cA=\setsuch{(x,t)\in\R\times[t_0,t_1]}{g_1(t)<x<g_2(t)}$
with continuous functions $g_1<g_2$.  Note that in this case,
$\tau_\cA$ is a stopping time\footnote{%
For a general Borel-measurable set $\cA$, the first exit time
$\tau_\cA$ is still a stopping time with respect to the canonical
filtration, completed by the null sets.}  
with respect to the canonical filtration
of $(\Omega, \cF, \fP)$ generated by $\{x_t\}_{t\geqs t_0}$.

Before turning to the precise statements of our results, let us
introduce some notations. We shall use
\begin{itemiz}
\item 
$\intpartplus{y}$ for $y\geqs 0$ to denote the smallest integer which
is greater than or equal to $y$, and 
\item 
$y\vee z$ and $y \wedge z$ to denote the maximum or minimum,
respectively, of two real numbers $y$ and $z$.
\item 
By $g(u)=\Order{u}$ we indicate that there exist $\delta>0$ and $K>0$
such that $g(u)\leqs K u$ for all $u\in[0,\delta]$, where $\delta$ and
$K$ of course do not depend on $\eps$ or $\sigma$. Similarly,
$g(u)=\orderone{}$ is to be understood as $\lim_{u\to0} g(u)=0$. From
time to time, we write $g(u)=\orderone{T}$ to indicate that choosing
{\it a priori} a sufficiently small $T$ allows to make the
corresponding term arbitrarily small for all $u$ from some
$T$-dependent interval. 
\end{itemiz}
Finally, let us point out that most estimates hold for small enough
$\eps$ only, and often only for $\fP$-almost all $\omega\in\Omega$. We
will stress these facts only when confusion might arise.  


\subsection{Nonbifurcating equilibria}
\label{ssec_rn}

We start by considering the nonlinear SDE \eqref{SDE}
in the case of $f$ admitting a nonbifurcating equilibrium branch. We will
assume that there exists an interval $I=[0,T]$ or $[0,\infty)$ such that the
following properties hold:
\begin{itemiz}
\item   there exists a function $x^\star:I\to \R$, called \defwd{equilibrium
curve}, such that
\begin{equation}
\label{rn1}
f(x^\star(t),t)=0
\qquad \forall t\in I;
\end{equation}

\item   $f$ is twice continuously differentiable with respect to $x$ and
$t$, with uniformly bounded derivatives, for all $t\in I$ and all $x$ in a
neighbourhood of $x^\star(t)$;

\item   the linearization of $f$ at $x^\star(t)$, defined as
\begin{equation}
\label{rn2}
a(t) = \sdpar fx(x^\star(t),t),
\end{equation}
is bounded away from zero, that is, there exists a constant $a_0>0$ such
that 
\begin{equation}
\label{rn3}
\abs{a(t)}\geqs a_0 \quad \forall t\in I.
\end{equation}
\end{itemiz}
In the deterministic case $\sigma=0$, the following result is
known~(see~\figref{fig1}a): 

\begin{theorem}[Deterministic case {\rm \cite{Tihonov,Grad}}]
\label{t_Tihonov}
Consider the equation
\begin{equation}
\label{ODEa}
\eps\dtot{x_t}t = f(x_t,t).
\end{equation}
There are constants $\eps_0, c_0, c_1>0$, depending only on $f$,
such that for $0<\eps\leqs\eps_0$,
\begin{itemiz}
\item   \eqref{ODEa} admits a particular solution $\xhatdet_t$ such that
\begin{equation}
\label{rn4b}
\abs{\xhatdet_t-x^\star(t)} \leqs c_1\eps
\quad \forall t\in I;
\end{equation}
\item   if $\abs{x_0 - x^\star(0)}\leqs c_0$ and $a(t)\leqs -a_0$ for all
$t\in I$ (that is, when $x^\star$ is a stable equilibrium), then the
solution $\xdet_t$ of \eqref{ODEa} with initial condition $\xdet_0=x_0$
satisfies
\begin{equation}
\label{rn4c}
\abs{\xdet_t - \xhatdet_t} \leqs \abs{x_0 - \xhatdet_0} \e^{-a_0 t/2\eps} 
\quad \forall t\in I.
\end{equation}
\end{itemiz}
\end{theorem}

\begin{remark}
\label{r_adiabatic}
The particular solution $\xhatdet$ is often called a \defwd{slow solution} or
\defwd{adiabatic solution} of equation \eqref{ODEa}. It is not unique in
general, as suggested by \eqref{rn4c}.
\end{remark}

We return now to the SDE \eqref{SDE} with $\sigma>0$. We need no additional
assumption on $\sigma$ in this section. However, the results are only
interesting when $\sigma = \orderone{\eps}$. 
Let us first consider the stable case, that is, we assume that $a(t)\leqs
-a_0<0$ for all $t\in I$. We assume that at $t=0$, $x_t$ starts at some
(deterministic) $x_0$ sufficiently close to $x^\star(0)$. Theorem
\ref{t_Tihonov} tells us that the deterministic solution $\xdet_t$ with
the same initial condition $\xdet_0=x_0$ reaches a neighbourhood of order
$\eps$ of $x^\star(t)$ exponentially fast.

We are interested in the stochastic process $y_t = x_t - \xdet_t$, which
describes the deviation due to noise from the deterministic solution
$\xdet$. It obeys the SDE
\begin{equation}
\label{rn5}
\6y_t = \frac1\eps \bigbrak{f(\xdet_t+y_t,t) - f(\xdet_t,t)} \6t 
+ \frac\sigma{\sqrt\eps} \6W_t, 
\qquad y_0 = 0.
\end{equation}
We will prove that $y_t$ remains in a neighbourhood of $0$ with high
probability. It is instructive to consider first the linearization of
\eqref{rn5} around $y=0$, which has the form
\begin{equation}
\label{rn6}
\6y^0_t = \frac1\eps \ba(t) y^0_t \6t + \frac\sigma{\sqrt\eps} \6W_t, 
\end{equation}
where 
\begin{equation}
\label{rn7}
\ba(t)  = \sdpar fx(\xdet_t,t) = a(t) + \Order{\eps} +
\bigOrder{\abs{x_0-x^\star(0)}\e^{-a_0t/2\eps}}.
\end{equation}
Taking $\eps$ and $\abs{x_0-x^\star(0)}$ sufficiently small, we may assume
the existence of constants $\ba_+\geqs \ba_->0$ such that $-\ba_+ \leqs
\ba(t) \leqs -\ba_-$ for all $t\in I$.
The solution of \eqref{rn6} with arbitrary initial condition $y^0_0$ is
given by
\begin{equation}
\label{rn8}
y^0_t = y^0_0 \e^{\balpha(t)/\eps} + 
\frac\sigma{\sqrt\eps} \int_0^t \e^{\balpha(t,s)/\eps}\6W_s, 
\qquad \balpha(t,s) = \int_s^t \ba(u)\6u,
\end{equation}
where we write $\balpha(t,0)=\balpha(t)$ for brevity. Note that
$\balpha(t,s)\leqs -\ba_-(t-s)$ whenever $t\geqs s$.  If 
$y^0_0$ has variance $v_0\geqs0$, then $y^0_t$ has variance
\begin{equation}
\label{rn9}
v(t) = v_0 \e^{2\balpha(t)/\eps} + \frac{\sigma^2}\eps \int_0^t
\e^{2\balpha(t,s)/\eps} \6s. 
\end{equation} 
Since the first term decreases exponentially fast, the initial variance
$v_0$ is \lq\lq forgotten\rq\rq\ as soon as $\e^{2\balpha(t)/\eps}$ is
small enough, which happens already for $t>\Order{\eps\abs{\log\eps}}$. 
For $y^0_0=0$, \eqref{rn8} implies in particular that for any $\delta>0$, 
\begin{equation}
\label{rn10}
\bigprobin{0,0}{\abs{y^0_t}\geqs\delta} \leqs \e^{-\delta^2/2v(t)}, 
\end{equation}
and thus the probability of finding $y^0_t$, at any given $t\in I$, outside
a strip of width much larger than $\sqrt{2v(t)}$ is very small. 

Our first main result states that the {\em whole path} $\set{x_s}_{0\leqs
s\leqs t}$ of the solution of the {\em nonlinear} equation \eqref{SDE} lies
in a similar strip with high probability. We only need to make one
concession: the width of the strip has to be bounded away from zero. 
Therefore, we define the strip as
\begin{equation}
\label{rn11}
\Bs(h) = \bigsetsuch{(x,t)\in\R\times I}{\abs{x-\xdet_t}<h\sqrt{\z(t)}},
\end{equation}
where 
\begin{equation}
\label{rn12}
\z(t) = \frac1{2\abs{\ba(0)}} \e^{2\balpha(t)/\eps} 
+ \frac1\eps \int_0^t \e^{2\balpha(t,s)/\eps} \6s.
\end{equation}
$\sigma^2\z$ can be interpreted as the variance \eqref{rn9} of the process
\eqref{rn8} starting with initial variance $v_0=\sigma^2/(2\abs{\ba(0)})$. 
We shall show in Lemma~\ref{l_ns1} that 
\begin{equation}
\label{rn13}
\z(t) = \frac1{2\abs{a(t)}} + \Order\eps +
\bigOrder{\abs{x_0-x^\star(0)}\e^{-a_0t/2\eps}}.
\end{equation}
Let $\tau_{\Bs(h)}$ denote the first exit time of $x_t$ from
$\Bs(h)$. 

\begin{theorem}[Stable case]
\label{t_stable}
There exist $\eps_0$, $d_0$ and $h_0$, depending only on $f$, such that for
$0<\eps\leqs\eps_0$, $h\leqs h_0$ and $\abs{x_0-x^\star(0)}\leqs d_0$, 
\begin{equation}
\label{rn14a}
\bigprobin{0,x_0}{\tau_{\Bs(h)}<t} \leqs C(t,\eps)
\exp\Bigset{-\frac12\frac{h^2}{\sigma^2}\bigbrak{1-\Order{\eps}-\Order{h}}},
\end{equation}
where
\begin{equation}
\label{rn14b}
C(t,\eps) = \frac{\abs{\balpha(t)}}{\eps^2} + 2.
\end{equation}
\end{theorem}

The proof, given in Section \ref{ssec_ns}, is divided into two main steps.
First, we show that an estimate of the form \eqref{rn14a}, but without the
term $\Order{h}$, holds for the solution of the linear equation
\eqref{rn6}. Then we show that whenever $\abs{y^0_s}<h\sqrt{\z(s)}$ for
$0\leqs s\leqs t$, one almost surely also has
$\abs{y_s}<h(1+\Order{h})\sqrt{\z(s)}$ for $0\leqs s\leqs t$.

\begin{remark}
\label{r_zeta}
The result of the preceding theorem remains true when
$1/2\abs{\ba(0)}$ in the definition \eqref{rn12} of $\z(t)$ is
replaced be an arbitrary $\z_0$, {\em provided} $\z_0>0$. The terms
$\Order{\cdot}$ may then depend on $\z_0$. Note that $\z(t)$ and
$\sigma^2 v(t)$ are both solutions of the same differential equation
$\eps z'=2\ba(t)z+1$, with possibly different initial conditions. If
$x_0-x^\star(0)=\Order{\eps}$, $\z(t)$ is an adiabatic solution (in
the sense of Theorem~\ref{t_Tihonov}) of the differential equation,
staying close to the equilibrium branch $z^\star=1/\abs{2\ba(t)}$. 
\end{remark}

The estimate \eqref{rn14a} has been designed for situations where
$\sigma\ll1$, and is useful for $\sigma\ll h\ll 1$. We expect the
exponent to be optimal in this case, but did not attempt to optimize the
prefactor $C(t,\eps)$, which leads to subexponential corrections. 
If we assume, for instance, that $\sigma=\eps^q$, $q>0$, and take $h=\eps^p$
with $0<p<q$, \eqref{rn14a} can be written as
\begin{equation}
\label{rn15}
\bigprobin{0,x_0}{\tau_{\Bs(h)}<t} \leqs 
(t+\eps^2) \exp\Bigset{-\frac1{2\eps^{2(q-p)}} \bigbrak{1-\Order{\eps} -
\Order{\eps^p} - \Order{\eps^{2(q-p)}\abs{\log\eps}}}}. 
\end{equation} 
The $t$-dependence of the prefactor is to be expected. It is due to the
fact that as time increases, the probability of $x_t$ escaping from a
neighbourhood of $\xdet_t$ also increases, but very slowly if $\sigma$
is small. The estimate \eqref{rn14a} shows that for a fraction $\gamma$ of
trajectories to leave the strip $\Bs(h)$, we have to wait {\em at least} for
a time $t_\gamma$ given by
\begin{equation}
\label{rn16}
\abs{\balpha(t_\gamma)} 
= \gamma\eps^2
\exp\Bigset{\frac12\frac{h^2}{\sigma^2}\bigbrak{1-\Order{\eps}-\Order{h}}} 
- 2\eps^2,
\end{equation}
which is compatible with results on the autonomous case.

Let us now consider the unstable case, that is, we now assume that the
linearization $a(t)=\sdpar fx(x^\star(t),t)$ satisfies $a(t)\geqs a_0>0$
for all $t\in I$. Theorem \ref{t_Tihonov} shows the existence of a
particular solution $\xhatdet_t$ of the deterministic equation \eqref{ODEa}
such that $\abs{\xhatdet_t-x^\star(t)} \leqs c_1\eps$ for all $t\in I$.  We
define $\ba(t)=\sdpar fx(\xhatdet_t,t) = a(t)+\Order{\eps}>0$ and
$\balpha(t) = \int_0^t \ba(s)\6s$.  

The linearization of \eqref{SDE} around $\xhatdet_t$ again admits a
solution of the form~\eqref{rn8}. In this case, however, the variance
\eqref{rn9} grows exponentially fast, and thus one expects the
probability of $x_t$ remaining close to $\xhatdet_t$ to be small. This
is the contents of the second main result of this section. We
introduce the set 
\begin{equation}
\label{rn20}
\Bu(h) = \biggsetsuch{(x,t)\in\R\times I}
{\abs{x-\xhatdet_t} < \frac h{\sqrt{2\ba(t)}}}
\end{equation}
and the first exit time $\tau_{\Bu(h)}$ of $x_t$ from $\Bu(h)$.

\begin{theorem}[Unstable case]
\label{t_unstable}
There exist $\eps_0$ and $h_0$, depending only on $f$, such that for all
$h\leqs \sigma\wedge h_0$, all $\eps\leqs\eps_0$ and all $x_0$ satisfying
$(x_0,0)\in\Bu(h)$, we have
\begin{equation}
\label{rn21}
\bigprobin{0,x_0}{\tau_{\Bu(h)}\geqs t} 
\leqs \sqrt{\e} 
\exp\Bigset{-\kappa \frac{\sigma^2}{h^2} \frac{\balpha(t)}\eps},
\end{equation} 
where $\kappa = \frac\pi{2{\e}} \bigpar{1-\Order{h}-\Order{\eps}}$.
\end{theorem}

The proof, given in Section~\ref{ssec_nu}, is based on a partition of
the interval $[0,t]$ into small intervals, and a comparison of the
nonlinear equation with its linearization on each interval.  

This result shows that $x_t$ is unlikely to remain in $\Bu(h)$ as
soon as $t\gg\eps\sigma^2/h^2$.  A major limitation of \eqref{rn21} is that
it requires $h\leqs\sigma$. Obtaining an estimate for larger $h$ is
possible, but requires considerably more work. We will provide such an
estimate in the more difficult, but also more interesting case of the
pitchfork bifurcation, see Theorem~\ref{t_escape} below.


\subsection{Pitchfork bifurcation}
\label{ssec_rb}

We now consider the SDE \eqref{SDE} in the case of $f$ undergoing a
pitchfork bifurcation. We will assume that
\begin{itemiz}
\item   $f$ is three times continuously differentiable with respect to $x$
and $t$ in a neighbourhood $\cN_0$ of $(0,0)$;
\item   $f(x,t) = -f(-x,t)$ for all $(x,t)\in\cN_0$;
\item   $f$ exhibits a supercritical pitchfork bifurcation at the origin,
i.e.
\begin{equation}
\label{rb1}
\sdpar{f}{x}(0,0) = 0, \qquad
\sdpar{f}{tx}(0,0) > 0  \qquad \text{and} \qquad 
\sdpar{f}{xxx}(0,0) < 0.
\end{equation} 
\end{itemiz}

The assumption that $f$ be odd is not necessary for the existence of a
pitchfork bifurcation. However, the deterministic system behaves very
differently if $x=0$ is not always an equilibrium. The most natural
situation in which $f(0,t)=0$ for all $t$ is the one where $f$ is odd.

By rescaling $x$ and $t$, we may arrange that $\sdpar f{tx}(0,0)=1$ and
$\sdpar f{xxx}(0,0) = -6$ as in the standard case $f(x,t)=tx-x^3$.  This
implies in particular that the linearization of $f$ at $x=0$ satisfies
\begin{equation}
\label{rb2}
a(t) = \sdpar fx(0,t) = t + \Order{t^2}.
\end{equation}
A standard result of bifurcation theory \cite{GH,IJ}
states that under these assumptions, there is a neighbourhood
$\cN\subset\cN_0$ of $(0,0)$ in which the only solutions of $f(x,t)=0$ are
the line $x=0$ and the curves
\begin{equation}
\label{rb3}
x = \pm x^\star(t), \qquad
x^\star(t) = \sqrt t \bigbrak{1+\orderone{t}}, \qquad
t\geqs 0.
\end{equation}
If $\cN$ is small enough,  the equilibrium $x=0$ is stable for $t<0$ and
unstable for $t>0$, while  $x = \pm x^\star(t)$ are stable equilibria with
linearization
\begin{equation}
\label{rb4}
a^\star(t) = \sdpar fx(x^\star(t),t) = -2t \bigbrak{1+\orderone{t}}.
\end{equation}
The only solutions of $\sdpar fx(x,t)=0$ in $\cN$ are the curves
\begin{equation}
\label{rb5}
x = \pm \bx(t), \qquad
\bx(t) = \sqrt {t/3} \bigbrak{1+\orderone{t}}, \qquad
t\geqs 0.
\end{equation}
If $f$ is four times continuously differentiable, the terms $\orderone{t}$ in
the last three equations can be replaced by $\Order{t}$. 

We briefly state what is known for the deterministic equation
\begin{equation}
\label{ODEb}
\eps\dtot {x_t}t = f(x_t,t),
\end{equation}
where we take an initial condition $(x_0,t_0)\in\cN$ with $x_0>0$ and
$t_0<0$, see~\figref{fig1}b. Observe that $\alpha(t,t_0) =
\int_{t_0}^t a(s) \6s$ is decreasing for $t_0<t<0$ and increasing for $t>0$.

\begin{definition}
\label{d_bifdelay}
The \defwd{bifurcation delay} is defined as
\begin{equation}
\label{bifdelay}
\Pi(t_0) = \inf\bigsetsuch{t>0}{\alpha(t,t_0)>0}, 
\end{equation}
with the convention $\Pi(t_0)=\infty$ if $\alpha(t,t_0)<0$ for all
$t>0$, for which $\alpha(t,t_0)$ is defined.
\end{definition}

One easily shows that $\Pi(t_0)$ is differentiable for $t_0$ sufficiently
close to $0$, and satisfies $\lim_{t_0\to 0-} \Pi(t_0)=0$ and $\lim_{t_0\to
0-} \Pi'(t_0)=-1$.

\begin{theorem}[Deterministic case]
\label{t_bifdelay}
Let $\xdet_t$ be the solution of \eqref{ODEb} with initial condition
$\xdet_{t_0}=x_0$. Then there exist constants $\eps_0$, $c_0$, $c_1$ depending
only on $f$, and times
\begin{equation}
\label{rb6}
\begin{split}
t_1 &= t_0 + \Order{\eps\abs{\log\eps}} \\
t_2 &= \Pi(t_1) = \Pi(t_0) - \Order{\eps\abs{\log\eps}} \\
t_3 &= \Pi(t_0) + \Order{\eps\abs{\log\eps}}
\end{split}
\end{equation}
such that, if $0<x_0\leqs c_0$, $0<\eps\leqs\eps_0$ and $(\xdet_t,t)\in\cN$, 
\begin{equation}
\label{rb7}
\begin{cases}
0 < \xdet_t \leqs c_1\eps \e^{\alpha(t,t_1)/\eps}
&\text{for $t_1 \leqs t \leqs t_2$} \\
\abs{\xdet_t-x^\star(t)} \leqs c_1\eps 
&\text{for $t \geqs t_3$.}
\end{cases}
\end{equation}
\end{theorem}

The proof is a straightforward consequence of differential inequalities,
see for instance \cite[Propositions 4.6 and 4.8]{Thesis}. 

We now consider the SDE \eqref{SDE} for $\sigma>0$. 
The results in this section are only interesting for
$\sigma=\order{\sqrt\eps}$, while one of them (Theorem~\ref{t_escape}) {\em
requires} a condition of the form $\sigma\abs{\log\sigma}^{3/2} =
\Order{\sqrt\eps}$ (where we have not tried to optimize the exponent $3/2$). 

Let us fix an initial
condition $(x_{t_0},t_0)\in\cN$ with $t_0<0$. For any $T\in(0,\abs{t_0})$,
we can apply Theorem \ref{t_stable} on the interval $[t_0,-T]$ to show that
$\abs{x_{-T}}$ is likely to be of order
$\sigma^{1-\delta} + c_1\eps\e^{\alpha(-T,t_1)/\eps}$ for any
$\delta>0$. We can also apply the theorem for $t>T$ to show that
the curves $\pm x^\star(t)$ attract nearby trajectories. Hence there
is no limitation in considering the SDE \eqref{SDE} in a domain of the
form $\abs{x}\leqs d$, $\abs{t}\leqs T$ where $d$ and $T$ can be taken
small (independently of $\eps$ and $\sigma$ of course!), with an
initial condition $x_{-T}=x_0$ satisfying $\abs{x_0}\leqs d$. 

We first show that $x_t$ is likely to remain small for $-T\leqs
t\leqs\sqrt\eps$. Actually, it turns out to be convenient to show that $x_t$
remains close to the solution $x_0\e^{\alpha(t,-T)/\eps}$ of the
linearization of \eqref{ODEb}. We define the ``variance-like'' function
\begin{equation}
\label{rb8}
\z(t) = \frac1{2\abs{a(-T)}} \e^{2\alpha(t,-T)/\eps} + 
\frac1\eps \int_{-T}^t \e^{2\alpha(t,s)/\eps} \6s.
\end{equation}
We shall show in Lemma~\ref{l_ps1} that for sufficiently small $\eps$,
there exist constants $c_{\pm}$ such that
\begin{align}
\label{rb9a}
\frac{c_-}{\abs t} &\leqs \z(t) \leqs \frac{c_+}{\abs{t}} 
&&\text{for $-T\leqs t\leqs -\sqrt\eps$,}\\
\label{rb9b}
\frac{c_-}{\sqrt\eps} &\leqs \z(t) \leqs \frac{c_+}{\sqrt\eps} 
&&\text{for $-\sqrt\eps\leqs t\leqs \sqrt\eps$.}
\end{align}
The function $\z(t)$ is used to define the strip
\begin{equation}
\label{rb10}
\cB(h) = \bigsetsuch{(x,t)\in[-d,d\mskip2mu]\times [-T,\sqrt\eps\,]}
{\abs{x-x_0\e^{\alpha(t,-T)/\eps}} < h\sqrt{\z(t)}}.
\end{equation}
Let $\tau_{\cB(h)}$ denote the first exit time of $x_t$ from $\cB(h)$.

\begin{theorem}[Behaviour for $t\leqs\sqrt\eps\mskip3mu$]
\label{t_before}
There exist constants $\eps_0$ and $h_0$, depending only on $f$, $T$
and $d$, such that for $0<\eps\leqs\eps_0$, $h\leqs h_0\sqrt\eps$,
$\abs{x_0}\leqs h/\eps^{1/4}$ and $-T\leqs t\leqs\sqrt\eps$,
\begin{equation}
\label{rb11a}
\bigprobin{-T,x_0}{\tau_{\cB(h)} < t}
\leqs C(t,\eps) \exp\biggset{-\frac12 \frac{h^2}{\sigma^2}
\biggbrak{1-r(\eps)-\biggOrder{\frac{h^2}\eps}}}
\end{equation}
where
\begin{equation}
\label{rb11b}
C(t,\eps) = \frac{\abs{\alpha(t,-T)} + \Order{\eps}}{\eps^2},
\end{equation}
and with $r(\eps)=\Order{\eps}$ for $-T\leqs t\leqs-\sqrt\eps$, 
and $r(\eps)=\Order{\sqrt\eps}$ for $-\sqrt\eps\leqs t\leqs\sqrt\eps$. 
\end{theorem}

The proof (given in Section \ref{ssec_psmallt}) and the interpretation of
this result are very close in spirit to those of Theorem \ref{t_stable}.
The only difference lies in the kind of $\eps$-dependence of
the error terms. The estimate \eqref{rb11a} is useful when $\sigma \ll
h \ll \sqrt\eps$, and shows that the typical spreading of paths around
the deterministic solution will slowly grow until $t=\sqrt\eps$, where
it is of order $\sigma/\eps^{1/4}$, see~\figref{fig2}.

Let us now examine what happens for $t\geqs\sqrt\eps$. We first show that
$x_t$ is likely to leave quite soon a suitably defined region $\cD$
containing the line $x=0$. The boundary of $\cD$ is defined through a
function $\tx(t)$, which can be chosen somewhat arbitrarily, but should lie
between $\bx(t)$ and $x^\star(t)$, in order to simplify the analysis of the
dynamics after $x_t$ has left $\cD$. A convenient definition is
\begin{equation}
\label{rb12}
\tx(t) = \sqrt\lambda \, x^\star(t),
\end{equation} 
where $\lambda$ is a free parameter. We need to assume, however, that
$\lambda\in(\frac13,\frac12)$. We now define
\begin{equation}
\label{rb13}
\cD = \bigsetsuch{(x,t)\in[-d,d\mskip2mu]\times[\sqrt\eps,T]}
{\abs x < \tx(t)}.
\end{equation}
Note that $\cD$ has the property that for all $(x,t)\in\cD$ with $x\neq
0$, 
\begin{equation}
\label{rb14}
\frac1x f(x,t) \geqs \kappa a(t)
\qquad
\text{with $\kappa = 1 - \lambda - \orderone{T}$.}
\end{equation}
Let $\tau_{\cD}$ denote the first exit time of $x_t$ from $\cD$.

\begin{theorem}[Escape from $\cD$]
\label{t_escape}
Let $(x_0,t_0)\in\cD$ and assume that
$\sigma\abs{\log\sigma}^{3/2}=\Order{\sqrt\eps}$. Then for $t_0\leqs
t\leqs T$,
\begin{equation}
\label{rb15}
\bigprobin{t_0,x_0}{\tau_{\cD}\geqs t}\leqs 
C_0 \,\tx(t)\sqrt{a(t)}\mskip2mu \frac{\abs{\log\sigma}}{\sigma} 
\biggpar{1+\frac{\alpha(t,t_0)}\eps}
\frac{\e^{-\kappa \alpha(t,t_0)/\eps}}{\sqrt{1-\e^{-2\kappa\alpha(t,t_0)/\eps}}},
\end{equation}
where $C_0>0$ is a (numerical) constant.
\end{theorem}

The proof of this result (given in Section \ref{ssec_pesc}) is by far the
most involved of the present work. We start by estimating, in a
similar way as in Theorem \ref{t_unstable}, the first exit time from a
strip $\cS$ of width slightly larger than $\sigma/\sqrt{a(s)}$. The
probability of returning to zero after leaving $\cS$ can be estimated;
it is small but not exponentially small. However, the probability of
neither leaving $\cD$ nor returning to zero is exponentially
small. This fact can be used to devise an iterative scheme that leads
to the exponential estimate \eqref{rb15}. 

We point out that for every subset $\cD'\subset\cD$, we have 
$\probin{t_0,x_0}{\tau_{\cD'}\geqs t} \leqs 
\probin{t_0,x_0}{\tau_{\cD}\geqs t}$, and thus \eqref{rb15} still
provides an upper bound for the first exit time from smaller sets.

Let us finally consider what happens after the path has left $\cD$ at
time $\tau = \tau_\cD$. One can deduce from the
definition \eqref{rb12} of $\tx(t)$ that for $\sqrt\eps\leqs t\leqs T$
and $\abs x\geqs \tx(t)$, 
\begin{equation}
\label{rb16}
\sdpar fx(x,t) \leqs \ta(t) = \sdpar fx(\tx(t),t) \leqs -\mm a(t)
\qquad
\text{with $\mm = 3\lambda - 1 - \orderone{T}$.}
\end{equation}
Let $\xdettau t$ denote the solution of the deterministic equation
\eqref{ODEb} starting in $\tx(t)$ at time $\tau$ (the case where one starts
at $-\tx(t)$ is obtained by symmetry). We shall show in
Proposition~\ref{p_app1} that $\xdettau t$ always remains between $\tx(t)$
and $x^\star(t)$, and approaches $x^\star(t)$ according to 
\begin{equation}
\label{rb17}
\xdettau t = x^\star(t) - \BigOrder{\frac\eps{t^{3/2}}} -
\bigOrder{\sqrt\tau \e^{-\mm\alpha(t,\tau)/\eps}}.
\end{equation}
Moreover, deterministic solutions starting at different times approach each
other like 
\begin{equation}
\label{rb18}
0 \leqs \xdetin{\sqrt\eps}t - \xdettau t \leqs 
\bigpar{\xdetin{\sqrt\eps}\tau - \tx(\tau)} \e^{-\mm\alpha(t,\tau)/\eps}
\qquad \forall t\in[\tau,T].
\end{equation}
The linearization of $f$ at $\xdettau t$ satisfies
\begin{equation}
\label{rb19}
\atau(t) = \sdpar fx(\xdettau t,t) 
= a^\star(t) + \BigOrder{\frac\eps t} +
\bigOrder{t\e^{-\mm\alpha(t,\tau)/\eps}}.
\end{equation}
For given $\tau$, we construct a strip $\cA^\tau(h)$ around $\xdettau{}$
of the form
\begin{equation}
\label{rb20}
\cA^\tau(h) = \bigsetsuch{(x,t)}
{\tau\leqs t\leqs T, \abs{x-\xdettau t} < h\sqrt{\ztau(t)}},
\end{equation}
where the function $\ztau(t)$ is defined by
\begin{equation}
\label{rb21}
\ztau(t) = \frac1{2\abs{\ta(\tau)}} \e^{2\alphatau(t,\tau)/\eps} 
+ \frac1\eps \int_\tau^t \e^{2\alphatau(t,s)/\eps}\6s, 
\qquad \alphatau(t,s) = \int_s^t \atau(u)\6u,
\end{equation}
and satisfies
\begin{equation}
\label{rb22}
\ztau(t) = \frac1{2\abs{a^\star(t)}} + \BigOrder{\frac\eps{t^3}} +
\BigOrder{\frac1t \e^{-\mm\alpha(t,\tau)/\eps}},
\end{equation}
cf.~Lemma~\ref{p_app2}. Let $\tau_{\cA^\tau(h)}$ denote the
first exit time of $x_t$ from $\cA^\tau(h)$.

\begin{theorem}[Approach to $x^\star$]
\label{t_approach}
There exist constants $\eps_0$ and $h_0$, depending only on $f$, $T$
and $d$, such that for $0<\eps\leqs\eps_0$, $h< h_0\tau$ and
$\tau\leqs t\leqs T$, 
\begin{equation}
\label{rb23a}
\bigprobin{\tau,\tx(\tau)}{\tau_{\cA^\tau(h)} < t}
\leqs C^\tau(t,\eps) \exp\Bigset{-\frac12 \frac{h^2}{\sigma^2}
\Bigbrak{1-\Order{\eps}-\BigOrder{\frac h\tau}}}
\end{equation}
where
\begin{equation}
\label{rb23b}
C^\tau(t,\eps) = \frac{\abs{\alphatau(t,\tau)}}{\eps^2} + 2
\leqs \frac1{\eps^2}\biggabs{\int_{\sqrt\eps}^t a^\star(s)\6s}+2.
\end{equation}
\end{theorem}

The proof is given in Section \ref{ssec_papp}. 
This result is useful for $\sigma\ll h\ll\tau$, and shows that the typical
spreading of paths around $\xdettau t$ is of order $\sigma/\sqrt{t}$,
see~\figref{fig2}.   


\subsection{Discussion}
\label{ssec_rd}

Let us now examine some of the consequences of these results.  First of
all, they allow to characterize the influence of additive noise on the
bifurcation delay. In the deterministic case, this delay is defined as the
first exit time from a strip of width $\eps$ around $x=0$, see Theorem
\ref{t_bifdelay}. A possible definition of the delay in the stochastic case
is thus the first exit time $\taudelay$ from a similar strip.  An
appropriate  choice for the width of the strip is $\tx(\sqrt\eps) =
\Order{\eps^{1/4}}$, since such a strip will contain $\cB(h)$ for every
admissible $h$, and the part of the strip with $t\geqs\sqrt\eps$ will be
contained in $\cD$.  Theorems \ref{t_before} and \ref{t_escape} then imply
that if $t\geqs\sqrt\eps$, 
\begin{align}
\label{rd1a}
\bigprobin{-T,x_0}{\taudelay < \sqrt\eps}
&\leqs C(\sqrt\eps,\eps) \e^{-\Order{\eps/\sigma^2}} \\
\label{rd2b}
\bigprobin{-T,x_0}{\taudelay \geqs t}
&\leqs C_0 \,\tx(t)\sqrt{a(t)}\mskip2mu \frac{\abs{\log\sigma}}{\sigma} 
\biggpar{1+\frac{\alpha(t,\sqrt\eps)}\eps}
\frac{\e^{-\kappa \alpha(t,\sqrt\eps)/\eps}}
{\sqrt{1-\e^{-2\kappa\alpha(t,\sqrt\eps)/\eps}}}.
\end{align}
If we choose $t$ in such a way that $\alpha(t,\sqrt\eps) = c\mskip1.5mu \eps
\abs{\log\sigma}$ for some $c>0$, the last expression reduces to
\begin{equation}
\label{rd3}
\bigprobin{-T,x_0}{\taudelay \geqs t} = \bigOrder{\sigma^{\kappa
c-1}\abs{\log\sigma}^2},
\end{equation}
which becomes small as soon as $c > 1/\kappa$. The bifurcation delay will
thus lie with overwhelming probability in the interval 
\begin{equation}
\label{rd4}
\bigbrak{\sqrt\eps, \bigOrder{\sqrt{\eps\abs{\log\sigma}}\,}}.
\end{equation}
Theorem \ref{t_approach} implies that for times larger than
$\Order{\sqrt{\eps\abs{\log\sigma}}\,}$, the paths are unlikely to return to
zero in a time of order~$1$. The wildest behaviour of the paths is to be
expected in the interval~\eqref{rd4}, because a region of instability
is crossed, where $\sdpar fx>0$. 

Our results on the pitchfork bifurcation require $\sigma\ll\sqrt\eps$,
while the estimate \eqref{rd4} is useful as long as $\sigma$ is not
exponentially small. We can thus distinguish three regimes, depending
on the noise intensity: 
\begin{enum}
\item[$\bullet$]
$\sigma\geqs\sqrt\eps$: A modification of Theorem~\ref{t_before} shows
that for $t<-\sigma$, the typical spreading of paths is of order
$\sigma/\sqrt{\abs t}$. Near the bifurcation point, the process is
dominated by noise, because the drift term $f\sim-x^3$ is too weak to
counteract the diffusion. Depending on the global structure of $f$, an
appreciable fraction of the paths might escape quite early from a
neighbourhood of the bifurcation point. In that situation, the notion
of bifurcation delay becomes meaningless.

\item[$\bullet$]        
$\e^{-1/\eps^p}\leqs\sigma\ll\sqrt\eps$ for some $p<1$: The
bifurcation delay lies in the interval \eqref{rd4} with high
probability, where $\sqrt{\eps\abs{\log\sigma}}\leqs \eps^{(1-p)/2}$
is still \lq\lq microscopic\rq\rq.

\item[$\bullet$]        
$\sigma\leqs\e^{-K/\eps}$ for some $K>0$: The noise is so small that
the paths remain concentrated around the deterministic solution for a time
interval of order $1$. The typical spreading is of order
$\sigma\sqrt{\z(t)}$, which behaves like $\sigma
\e^{\alpha(t)/\eps}/\eps^{1/4}$ for $t\geqs\sqrt\eps$, see
Lemma~\ref{l_ps1}. Thus the paths remain close to the origin until
$\alpha(t) \simeq \eps\abs{\log\sigma}\geqs K$. If
$\eps\abs{\log\sigma}>\alpha(\Pi(t_0)) = \abs{\alpha(t_0)}$, they 
follow the deterministic solution which makes a quick transition to
$x^\star(t)$ at $t=\Pi(t_0)$.
\end{enum}
The expression \eqref{rd4} characterizing the delay is in accordance with
experimental results \cite{TM,SMC}, and with the approximate calculation of
the last crossing of zero \cite{JL}. The numerical results in \cite{Gaeta},
which are fitted, at $\eps=0.01$, to $\taudelay\simeq\sigma^{0.105}$ for
weak noise and $\taudelay\simeq\e^{-851\,\sigma}$ for strong noise, seem
rather mysterious. Finally, the results in~\cite{Kuske}, who
approximates the probability density by a Gaussian centered at the
deterministic solution, can obviously only apply to the regime of
exponentially small noise. 

Another interesting question is how fast the paths concentrate near
the equilibrium branches $\pm x^\star(t)$. The deterministic
solutions, starting at $\tx(t_0)$ at some time $t_0>0$, all track
$x^\star(t)$ at a distance which is asymptotically of order
$\eps/t^{3/2}$. Therefore, we can choose one of them, say
\smash{$\xdetin{\sqrt\eps}t$}, and measure the distance of  $x_t$ from
that deterministic solution. We restrict our attention to those paths
which are still in a neighbourhood of the origin at time $\sqrt\eps$,
as most paths are. We want to show that for suitably chosen
$t_1\in(\sqrt\eps,t)$ and $\Delta\in(0,t)$, most paths will leave
$\cD$ until time $t_1$ and reach a $\delta$-neighbourhood of
\smash{$\xdetin{\sqrt\eps}t$} at time $\tau_\cD+\Delta$. Let us
estimate
\begin{align}
\label{rd5}
&\biggprobin{\sqrt\eps,x_{\sqrt\eps}}
{\biggpar{\tau_\cD<t_1,\ \sup_{s\in[\tau_\cD+\Delta,t]}
\Bigabs{\abs{x_s}-\xdetin{\sqrt\eps}s}<\delta}^{\math c}} \\  
\nonumber
&\quad
\leqs \bigprobin{\sqrt\eps,x_{\sqrt\eps}}{\tau_\cD\geqs t_1} 
+ \biggexpecin{\sqrt\eps,x_{\sqrt\eps}}{\indexfct{\tau_\cD<t_1}\;
\biggprobin{\tau_\cD,\tx(\tau_\cD)}{\sup_{s\in[\tau_\cD+\Delta,t]}
\abs{x_s-\xdetin{\sqrt\eps}s}\geqs\delta}}.
\end{align}
The first term decreases roughly like
$\sigma^{-1} \e^{-\kappa\alpha(t_1,\sqrt\eps)/\eps}$ and becomes small
as soon as $\alpha(t_1,\sqrt\eps)\gg\eps\abs{\log\sigma}$. The second
summand is bounded above by 
\begin{equation}
\label{rd6}
\text{\it const\/}\;
\Bigexpecin{\sqrt\eps,x_{\sqrt\eps}}{\indexfct{\tau_\cD<t_1}\;
\exp\Bigset{-\frac{t^2}{\sigma^2} 
\Bigbrak{\delta-\bigOrder{\sqrt{\tau_\cD}
\e^{-\y\alpha(\tau_\cD+\Delta,\tau_\cD)/\eps}}^2}}}.
\end{equation}
Therefore, $\delta$ should be large compared to $\sigma/t$ and we
also need that $\Delta$ is at least of order
$\Order{\sqrt{\eps\abs{\log\sigma}}}$. Then we see that after a time of order
$\Order{\sqrt{\eps\abs{\log\sigma}}}$, the typical paths will have left
$\cD$ and, after another time of the same order, will reach a
neighbourhood of \smash{$\xdetin{\sqrt\eps}t$}, which scales with
$\sigma/t$.

Finally, we can also estimate the probability of reaching the positive
rather than the negative branch. Consider $x_s$, starting in $x_0$ at
time $t_0<0$, and let $t>0$. Without loss of generality, we may assume that
$x_0>0$. The symmetry of $f$ implies 
\begin{equation}
\label{rd7}
\bigprobin{t_0,x_0}{x_t\geqs 0}
= 1 - \frac12 \bigprobin{t_0,x_0}{\exists\, s\in[t_0,t): x_s=0},
\end{equation}
and therefore it is sufficient to estimate the probability for $x_s$
to reach zero before time zero, for instance. We linearize the
SDE~\eqref{SDE} and use the fact that the solution $x^0_s$ of the
linearized equation 
\begin{equation}
\label{rd8}
\6x^0_s=\frac1\eps a(s)x^0_s \6s +\frac\sigma{\sqrt\eps} \6W_s,
\qquad 
x^0_{t_0}=x_0
\end{equation}
satisfies $x_s\leqs x^0_s$ as long as $x_s$ does not reach zero. For
the Gaussian process $x^0_s$ we know 
\begin{equation}
\label{rd9}
\bigprobin{t_0,x_0}{\exists\, s\in[t_0,t): x^0_s=0}
= 2\Bigpar{1-\bigprobin{t_0,x_0}{x^0_t\geqs 0}} 
= 1 - \frac1{\sqrt{2\pi}}\int_{-u(t)}^{u(t)} \e^{-y^2/2}\6y, 
\end{equation}
where $u(t)=x_0\e^{\alpha(t,t_0)/\eps}/\sqrt{v(t,t_0)}$ and $v(t,t_0)$
denotes the variance of $x^0_t$. For $t=0$, $u(0)$ is of order
$x_0\eps^{1/4}\sigma^{-1}\e^{-\text{\it const\/}\;t_0^2/\eps}$, see
Lemma~\ref{l_ps1}. Thus the probability in~\eqref{rd9} is
exponentially close to one for small $\eps$, and we conclude that the
probability for $x_t$ to reach the positive branch rather than the
negative one is exponentially close to $1/2$.


\section{The motion near nonbifurcating equilibria}
\label{sec_nobif}


In this section we consider the nonlinear SDE
\begin{equation}
\label{np1}
\6x_t = \frac{1}{\eps} f(x_t,t) \6t +
\frac{\sigma}{\sqrt{\eps}} \6W_t
\end{equation}
under the assumptions
\begin{itemiz}
\item   $t\in I = [0,T]$ or $[0,\infty)$;
\item   there exists an \defwd{equilibrium curve} $x^\star:I\to\R$ such that 
\begin{equation}
\label{np2}
f(x^\star(t),t) = 0
\quad \forall t\in I;
\end{equation}
\item   there is a constant $d>0$ such that $f$ is twice continuously
differentiable with respect to $x$ and $t$ for $\abs{x-x^\star(t)}\leqs d$
and $t\in I$, with $\abs{\sdpar{f}{xx}(x,t)}$ uniformly bounded by $2M>0$ in
that domain;
\item   there is a constant $a_0>0$ such that $a(t) =
\sdpar{f}{x}(x^\star(t),t)$ satisfies 
\begin{equation}
\label{np3}
\abs{a(t)} \geqs a_0
\quad \forall t\in I.
\end{equation}
\end{itemiz}
We do not need any assumptions on $\sigma>0$, but our results are of
interest only for $\sigma=\orderone{\eps}$.

In Section~\ref{ssec_ns} we consider the stable case, corresponding to
$a(t)\leqs -a_0<0$ for all $t\in I$. We first analyse the
linearization of~\eqref{np1} around a given deterministic
solution. Proposition~\ref{l_ns3} shows that the solutions of the
linearized equation are likely to remain in a strip of width
$h\sqrt{\z(t)}$ around the deterministic solution. Here $\z(t)$ is
related to the variance and will be analyzed in
Lemma~\ref{l_ns1}. Proposition~\ref{l_ns4} allows to 
compare the trajectories of the linear and the nonlinear equation, and
thus completes the proof of Theorem~\ref{t_stable}.

In Section~\ref{ssec_nu}, we consider the unstable case, i.\,e. $a(t)\geqs
a_0>0$ for all $t\in I$. Theorem~\ref{t_unstable} is equivalent to
Proposition~\ref{l_nu1}, which is again based on a comparison of
solutions of the nonlinear equation~\eqref{np1} and its linearization
around a given deterministic solution.


\subsection{Stable case}
\label{ssec_ns}

We first consider the case of a stable equilibrium, that is, we assume
that $a(t)\leqs -a_0$ for all $t\in I$. We will assume that the
stochastic process $x_t$, given by the SDE~\eqref{np1}, starts at time
$t=0$ in $x_0$. By Theorem~\ref{t_Tihonov}, there exists a $c_0>0$
such that the deterministic solution $\xdet$ of~\eqref{ODEa} with
initial condition $\xdet_0=x_0$ satisfies 
\begin{equation}
\label{ns0}
\abs{\xdet_t - x^\star(t)} \leqs 2c_1\eps + \abs{x_0-x^\star(0)}
\e^{-a_0t/2\eps}
\qquad \forall t\in I,
\end{equation}
provided $\abs{x_0-x^\star(0)}\leqs c_0$. We are interested in the
stochastic process $y_t = x_t - \xdet_t$, which describes the
deviation due to noise from the deterministic solution $\xdet$. It
obeys an SDE of the form 
\begin{equation}
\label{ns1}
\6y_t = \frac1\eps \bigbrak{\ba(t)y_t + \bb(y_t,t)} \6t +
\frac{\sigma}{\sqrt\eps} \6W_t,
\qquad y_0 = 0,
\end{equation}
where we have introduced the notations
\begin{equation}
\label{ns2}
\begin{split}
\ba(t) &= \ba_\eps(t) = \sdpar{f}{x}(\xdet_t,t) \\
\bb(y,t) &= \bb_\eps(y,t) = f(\xdet_t+y,t) - f(\xdet_t,t) - \ba(t)y.
\end{split}
\end{equation}
Taking $\eps$ and $\abs{x_0-x^\star(0)}$ sufficiently small, we may assume
that there exists a constant $\bd>0$ such that $\abs{\xdet_t+y-x^\star(t)}
\leqs d$ whenever $\abs{y}\leqs \bd$. It follows from Taylor's formula that
for all $(y,t)\in [-\bd,\bd\mskip2mu]\times I$, 
\begin{align}
\label{ns3a}
\abs{\bb(y,t)} &\leqs My^2 \\
\label{ns3b}
\abs{\ba(t)-a(t)} &\leqs M \bigpar{2c_1\eps +
\abs{x_0-x^\star(0)}\e^{-a_0t/2\eps}} 
\end{align}
By again taking $\eps$ and $\abs{x_0-x^\star(0)}$ sufficiently small,
we may further assume that there are constants $\ba_+\geqs\ba_->a_0/4$
such that
\begin{equation}
\label{ns3c}
-\ba_+ \leqs \ba(t) \leqs -\ba_-
\qquad \forall t\in I.
\end{equation}
Finally, the relation $\ba'(t) =
\sdpar{f}{xt}(\xdet_t,t) + \sdpar{f}{xx}(\xdet_t,t)\frac1\eps
f(\xdet_t,t)$ implies the existence of a constant $c_2>0$ such that
\begin{equation}
\label{ns3d}
\abs{\ba'(t)} \leqs c_2\Bigpar{1+\abs{x_0-x^\star(0)}
\frac{\e^{-a_0t/2\eps}}\eps}.
\end{equation}

Our analysis will be based on a comparison between solutions of \eqref{ns1}
and those of the linearized equation
\begin{equation}
\label{ns4}
\6y^0_t = \frac1\eps \ba(t) y^0_t \6t + \frac\sigma{\sqrt\eps} \6W_t, 
\qquad y^0_0 = 0.
\end{equation}
Its solution is given by 
\begin{equation}
\label{ns5}
y^0_t = \frac\sigma{\sqrt\eps} \int_0^t \e^{\balpha(t,s)/\eps}\6W_s, 
\qquad \balpha(t,s) = \int_s^t \ba(u)\6u.
\end{equation}
We will write $\balpha(t,0)=\balpha(t)$ for brevity. The Gaussian
random variable $y^0_t$ has mean zero and variance 
\begin{equation}
\label{ns6}
v(t) = \frac{\sigma^2}\eps \int_0^t \e^{2\balpha(t,s)/\eps} \6s.
\end{equation}
Note that \eqref{ns3c} implies that $\balpha(t,s)\leqs -\ba_-(t-s)$ whenever
$t\geqs s$, which implies in particular, that $v(t)$ is not larger than
$\sigma^2/2\ba_-$. We can, however, derive a more precise bound, which
is useful when $\eps$ and $\e^{-a_0 t/2\eps}$ are small. To do so, we
introduce the function
\begin{equation}
\label{ns6a}
\z(t) = \frac1{2\abs{\ba(0)}} \e^{2\balpha(t)/\eps} + \frac1\eps \int_0^t
\e^{2\balpha(t,s)/\eps}\6s,
\qquad \text{where\ } \balpha(t)=\balpha(t,0).
\end{equation}
Note that $v(t)\leqs \sigma^2 \z(t)$, and that both functions differ
by a term which becomes negligible as soon as
$t>\Order{\eps\abs{\log\eps}}$. The behaviour of $\z(t)$ is
characterized in the following lemma.

\begin{lemma}
\label{l_ns1}
The function $\z(t)$ satisfies the following relations for all $t\in I$. 
\begin{gather}
\label{ns7a}
\z(t) = \frac1{2\abs{\ba(t)}} + \Order\eps +
\bigOrder{\abs{x_0-x^\star(0)}\e^{-a_0t/2\eps}} \\
\label{ns7b}
\frac1{2\ba_+} \leqs \z(t) \leqs \frac1{2\ba_-} \\
\label{ns7c}
\z'(t) \leqs \frac1\eps
\end{gather}
\end{lemma}
\begin{proof}
By integration by parts, we obtain that
\begin{equation}
\label{ns8a}
\z(t)  
= \frac1{-2\ba(t)} - \frac12\int_0^t
\frac{\ba'(s)}{\ba(s)^2} 
\e^{2\balpha(t,s)/\eps} \6s.
\end{equation}
Using \eqref{ns3c} and \eqref{ns3d} we get
\begin{align}
\nonumber
&\Bigabs{\int_0^t\frac{\ba'(s)}{\ba(s)^2} \e^{2\balpha(t,s)/\eps} \6s}\\
\nonumber
& \qquad \leqs \frac{c_2}{\ba_-^2} \int_0^t \e^{-2\ba_-(t-s)/\eps}\6s 
+ \frac{c_2}{\ba_-^2} \frac{\abs{x_0-x^\star(0)}}\eps \int_0^t
\e^{[-2\ba_-(t-s)-a_0s/2]/\eps}\6s\\
&\qquad \leqs \frac{c_2}{2\ba_-^3}\eps + \frac{c_2}{\ba_-^2}
\frac{\abs{x_0-x^\star(0)}}{2\ba_- -a_0/2} \e^{-a_0t/2\eps},
\label{ns8b}
\end{align}
which proves \eqref{ns7a}.
We now observe that $\z(t)$ is a solution of the linear ODE
\begin{equation}
\label{ns8c}
\dtot\z t = \frac1\eps \bigpar{2\ba(t)\z + 1}, 
\qquad \z(0) = \frac1{2\abs{\ba(0)}}.
\end{equation}
Since $\z(t)>0$ and $\ba(t)<0$, we have $\z'(t)\leqs1/\eps$. 
We also see that $\z'(t)\geqs 0$ whenever $\z(t)\leqs 1/2\ba_+$ and 
$\z'(t)\leqs 0$ whenever $\z(t)\geqs 1/2\ba_-$. Since $\z(0)$ belongs to the
interval $[1/2\ba_+,1/2\ba_-]$, $\z(t)$ must remain in this interval for all
$t$.
\end{proof}

As we have already seen in~\eqref{rn10}, the probability of finding
$y^0_t$ outside a strip of width much larger than $\sqrt{2v(t)}$ is
very small. By Lemma~\ref{l_ns1}, we now know that $\sqrt{2v(t)}$
behaves approximately like $\sigma\abs{a(t)}^{-1/2}$. One of the key
points of the present work is to show that the {\em whole path}
$\set{y_s}_{0\leqs s\leqs t}$ remains in a strip of similar width
with high probability. The strip will be defined with the help of
$\zeta(t)$ instead of $v(t)$, because we need the width to be
bounded away from zero, even for small $t$.

To investigate $y^0_t$ we need to estimate the stochastic integral
from~\eqref{ns5}. Lemma~\ref{p_ns1} in the appendix provides the
estimate
\begin{equation}
\label{ns10z}
\Bigprob{\sup_{0\leqs s \leqs t} \int_0^s \varphi(u) \6W_u \geqs \delta} 
\leqs \exp\biggset{-\frac{\delta^2}{2\int_0^t \varphi(u)^2\6u}}
\end{equation}
for Borel-measurable deterministic functions $\varphi(u)$. 
Unfortunately, this estimate cannot be applied directly, because
in~\eqref{ns5}, the integrand depends explicitly on the upper
integration limit. This is why we introduce a partition of the
interval $[0,t]$. 

\begin{lemma}
\label{l_ns2}
Let $\rho:I\to\R_+$ be a measurable, strictly positive function. Fix
$K\in\N$, and let $0=u_0\leqs u_1 <\dots< u_K=t$ be a partition of
the interval $[0,t]$. Then
\begin{equation}
\label{ns11}
\Bigprobin{0,0}{\sup_{0\leqs s\leqs t} \frac{\abs{y^0_s}}{\rho(s)} \geqs h} 
\leqs 2\sum_{k=1}^K P_k,
\end{equation}
where
\begin{equation}
\label{ns11a}
P_k = \exp\biggset{-\frac12 \frac{h^2}{\sigma^2}  
\Bigpar{\inf_{u_{k-1}\leqs s\leqs u_k} \rho(s)^2 \e^{2\balpha(u_k,s)/\eps}}
\Bigpar{\frac1\eps \int_0^{u_k}\e^{2\balpha(u_k,s)/\eps}\6s}^{-1}}.
\end{equation}
\end{lemma}
\begin{proof}
We have 
\begin{align}
\label{ns12}
\Bigprobin{0,0}{&\sup_{0\leqs s\leqs t}
\frac{\abs{y^0_s}}{\rho(s)} \geqs h} \\
\nonumber
&= \Bigprobin{0,0}{\sup_{0\leqs s\leqs t} \frac1{\rho(s)} \Bigabs{\int_0^s
\e^{\balpha(s,u)/\eps}\6W_u} \geqs \frac{h\sqrt\eps}{\sigma}}\\
\nonumber
&= \Bigprobin{0,0}{\exists k \in\set{1,\dots,K}: 
\sup_{u_{k-1}\leqs s\leqs u_k} \frac1{\rho(s)} \Bigabs{\int_0^s
\e^{\balpha(s,u)/\eps} \6W_u} \geqs \frac{h\sqrt\eps}{\sigma}} \\
\nonumber
&\leqs 2\sum_{k=1}^{K} \Bigprobin{0,0}
{\sup_{u_{k-1}\leqs s\leqs u_k} \int_0^s
\e^{-\balpha(u)/\eps} \6W_u \geqs \frac{h\sqrt\eps}{\sigma} 
\inf_{u_{k-1}\leqs s\leqs u_k} \rho(s) \e^{-\balpha(s)/\eps}}.
\end{align}
Applying Lemma~\ref{p_ns1} to the last expression,
we obtain \eqref{ns11}.
\end{proof}

We are now ready to derive an upper bound for the probability that $y^0_s$
leaves a strip of appropriate width $h\rho(s)$ before time $t$. Taking 
$\rho(s)=\sqrt{\z(s)}$ will be a good choice since it leads to 
approximately constant $P_k$ in \eqref{ns11}. 

\begin{prop}
\label{l_ns3}
There exists an $r=r(\ba_+,\ba_-)$ such that 
\begin{equation}
\label{ns13}
\Bigprobin{0,0}{\sup_{0\leqs s\leqs t}\frac{\abs{y^0_s}}{\sqrt{\z(s)}}
\geqs h}  
\leqs C(t,\eps) \exp\Bigset{-\frac12\frac{h^2}{\sigma^2} (1-r\eps)},
\end{equation}
where
\begin{equation}
\label{ns14}
C(t,\eps) = \frac{\abs{\balpha(t)}}{\eps^2} + 2.
\end{equation}
\end{prop}
\begin{proof}
Let
\begin{equation}
\label{ns14a}
K = \biggintpartplus{\frac{\abs{\balpha(t)}}{2\eps^2}}.
\end{equation}
For $k = 1,\dots,K-1$, we define the partition times $u_k$ by the relation
\begin{equation}
\label{ns14b}
\abs{\balpha(u_k)} = 2\eps^2 k,
\end{equation}
which is possible since $\balpha(t)$ is continuous and decreasing. This
definition implies in particular that $\balpha(u_k,u_{k-1})=-2\eps^2$
and, therefore, $u_k-u_{k-1}\leqs 2\eps^2/\ba_-$. Bounding the
integral in \eqref{ns11a} by $\z(u_k)$, we obtain
\begin{equation}
\label{ns14c}
P_k \leqs \exp\Bigset{-\frac12 \frac{h^2}{\sigma^2} \inf_{u_{k-1}\leqs
s\leqs u_k} \frac{\z(s)}{\z(u_k)} \e^{2\balpha(u_k,s)/\eps}}.
\end{equation}
We have $\e^{2\balpha(u_k,s)/\eps}\geqs\e^{-4\eps}$ and 
\begin{equation}
\label{ns14d}
\z(s)-\z(u_k) = -\int_s^{u_k} \z'(u)\6u \geqs -\frac{u_k-s}\eps.
\end{equation}
Since $\z(u_k)\geqs 1/2\ba_+$, this implies
\begin{equation}
\label{ns14e}
P_k \leqs \exp\Bigset{-\frac12\frac{h^2}{\sigma^2}
\Bigpar{1-4\frac{\ba_+}{\ba_-}\eps}\e^{-4\eps}},
\end{equation}
and the result follows from Lemma \ref{l_ns2}.
\end{proof}

\begin{remark}
\label{r_ns2}
If we only assume that $\ba$ is Borel-measurable with $\ba(t)\leqs -\ba_-$ for
all $t\in I$, we still have 
\begin{equation}
\label{ns16}
\Bigprobin{0,0}{\sup_{0\leqs s\leqs t}\abs{y^0_s} \geqs h/\sqrt{2\ba_-}}
\leqs C(t,\eps) \exp\Bigset{-\frac12\frac{h^2}{\sigma^2} \e^{-4\eps}}.
\end{equation}
To prove this, we choose the same partition as before and bound the
integral in \eqref{ns11a} by $\eps/2\ba_-$.
\end{remark}

We now return to the nonlinear equation \eqref{ns1}, the solutions of which
we want to compare to those of its linearization \eqref{ns4}. To this
end, we introduce the events
\begin{align}
\label{ns17b}
\Omega_t(h) &= \Bigsetsuch{\w}{\bigabs{y_s(\w)} 
< h \sqrt{\z(s)} \; \forall s \in [0,t]}\\
\label{ns17a}
\Omega^0_t(h) &= \Bigsetsuch{\w}{\bigabs{y^0_s(\w)} 
< h \sqrt{\z(s)} \; \forall s \in [0,t]}.
\end{align}
Proposition~\ref{l_ns3} gives us an upper bound on the probability of the
complement of $\Omega^0_t(h)$. The key point to control the nonlinear case
is a relation between the sets $\Omega_t$ and $\Omega^0_t$ (for
slightly different values of $h$). This is done in
Proposition~\ref{l_ns4} below. 

\begin{notation}
\label{n_ns1}
For two events $\Omega_1$ and $\Omega_2$, we write $\Omega_1\subas\Omega_2$
if\/ $\fP$-almost all $\w\in\Omega_1$ belong to $\Omega_2$.
\end{notation}

\begin{prop}
\label{l_ns4}
Let $\gamma=2\sqrt{2\ba_+}\,M/\ba_-^2$ and assume that $h<\bd\sqrt{\ba_-/2}
\wedge \gamma^{-1}$. Then
\begin{align}
\label{ns18a}
\Omega_t(h) & \subas \Omega^0_t \Bigpar{\bigbrak{1+\frac\gamma4 h}h} \\
\Omega^0_t(h) & \subas \Omega_t \Bigpar{\bigbrak{1+\gamma h}h}.
\label{ns18b}
\end{align}
\end{prop}
\goodbreak
\begin{proof}\hfill
\begin{enum}
\item
The difference $z_s = y_s - y^0_s$ satisfies
\begin{equation}
\label{p:ns4:1}
\dtot{z_s}{s} = \frac{1}{\eps} \bigbrak{\ba(s)z_s + \bb(y^0_s+z_s,s)}
\end{equation}
with $z_0 = 0$ $\fP$-a.s.
Now,
\begin{equation}
\label{p:ns4:3neu}
z_s = \frac1\eps \int_0^s \e^{\balpha(s,u)/\eps} \bb(y^0_u+z_u,u) \6u,
\end{equation}
which implies
\begin{equation}
\label{p:ns4:4neu}
\abs{z_s} \leqs \frac{1}{\eps} 
\int_0^s\e^{\balpha(s,u)/\eps} \abs{\bb(y_u,u)} \6u
\end{equation}
for all $s\in[0,t]$. 
\item
Let us assume that $\w\in\Omega_t(h)$. Then we have for all $s\in[0,t]$ 
\begin{equation}
\label{p:ns4:6}
\abs{y_s(\w)} \leqs h\sqrt{\z(s)} \leqs
\frac{h}{\sqrt{2\ba_-}} \leqs \frac\bd2,
\end{equation}
and thus by \eqref{p:ns4:4neu},
\begin{equation}
\label{p:ns4:7neu}
\abs{z_s(\w)} \leqs \frac1\eps \int_0^s \e^{\balpha(s,u)/\eps}
\frac{Mh^2}{2\ba_-} \6u.
\end{equation}
The integral on the right-hand side can be estimated by~\eqref{ns7b},
yielding
\begin{equation}
\label{p:ns4:5}
\frac1\eps
\int_0^s \e^{\balpha(s,u)/\eps}\6u \leqs 2\z_{2\eps}(s) \leqs \frac1{\ba_-}. 
\end{equation}
Therefore,
\begin{equation}
\label{p:ns4:7neub}
\abs{z_s(\w)} \leqs \frac{Mh^2}{2\ba_-^2}
\leqs \frac{M\sqrt{\ba_+}\,h}{\sqrt2 \ba_-^2} h \sqrt{\z(s)},
\end{equation}
which proves \eqref{ns18a} because $\abs{y^0_s(\w)} \leqs
\abs{y_s(\w)} + \abs{z_s(\w)}$. 

\item
Let us now assume that $\w\in\Omega^0_t(h)$. Then we have 
$\abs{y^0_s(\w)} \leqs \bd/2$ for all $s\in[0,t]$ as in~\eqref{p:ns4:6}.
For $\delta=\gamma h$, we have $\delta<1$ by assumption. 
We consider the first exit time
\begin{equation}
\label{p:ns4:9}
\tau = \inf\bigsetsuch{s\in[0,t]}{\abs{z_s} \geqs \delta
h\sqrt{\z(s)}} \in [0,t]\cup\set{\infty}
\end{equation}
and the event
\begin{equation}
\label{p:ns4:10}
A = \Omega^0_t \cap \bigsetsuch{\w}{\tau(\w)<\infty}.
\end{equation}
If $\w\in A$, then for all $s\in[0,\tau(\w)]$, we have
$\abs{y_s(\w)}\leqs (1+\delta)h\sqrt{\z(s)}\leqs\bd$, and thus by
\eqref{p:ns4:4neu} and~\eqref{p:ns4:5},
\begin{equation}
\label{p:ns4:11}
\abs{z_s(\w)} \leqs \frac{1}{\eps} 
\int_0^s\e^{\balpha(s,u)/\eps} \frac{M(1+\delta)^2h^2}{2\ba_-}\6u 
\leqs \frac{M(1+\delta)^2h^2}{2\ba_-^2} 
< \delta h \sqrt{\z(s)}.
\end{equation}
However, by the definition of $\tau$, we have
$\abs{z_{\tau(\w)}(\w)}=\delta h\sqrt{\z(\tau(\w))}$,
which contradicts \eqref{p:ns4:11} for $s=\tau(\w)$. Therefore
$\prob{A}=0$, which implies that for almost all $\w\in \Omega^0_t$, we have
$\abs{z_s(\w)} < \delta h\sqrt{\z(s)}$ for all $s\in[0,t]$, and hence
\begin{equation}
\label{p:ns4:12}
\abs{y_s(\w)} < (1+\delta) h \sqrt{\z(s)} \quad \forall
s\in[0,t] 
\end{equation}
for these $\w$, which proves \eqref{ns18b}.
\qed
\end{enum}
\renewcommand{\qed}{}
\end{proof}

We close this subsection with a corollary which is
Theorem~\ref{t_stable}, restated in terms of the process $y_t$. 

\begin{cor}
\label{c_ns}
There exist $h_0$ and $\eps_0$, depending only on $f$, such that
for $\eps<\eps_0$ and $h<h_0$, 
\begin{equation}
\label{cns3}
\Bigprobin{0,0}{\sup_{0\leqs s\leqs t}
\frac{\abs{y_s}}{\sqrt{\z(s)}}>h} 
\leqs C(t,\eps)
\exp\Bigset{-\frac12\frac{h^2}{\sigma^2}\bigbrak{1-\Order{\eps}-\Order{h}}}.
\end{equation}
\end{cor}
\begin{proof}
By Proposition \ref{l_ns4} and Proposition~\ref{l_ns3}, 
\begin{equation}
\label{p:cns:3}
\begin{split}
\Bigprobin{0,0}{\sup_{0\leqs s\leqs t}
\frac{\abs{y_s}}{\sqrt{\z(s)}}>h} 
&\leqs
\Bigprobin{0,0}{\sup_{0\leqs s\leqs t}
\frac{\abs{y^0_s}}{\sqrt{\z(s)}}>h_1} \\
&\leqs C(t,\eps) \exp\Bigset{-\frac12 \frac{h_1^2}{\sigma^2}(1-r\eps)},
\end{split}
\end{equation}
where $h=(1+\gamma h_1)h_1$, which implies
\begin{equation}
\label{p:cns:4}
h_1 = \frac1{2\gamma}\bigbrak{\sqrt{1+4\gamma h}-1} \geqs h[1-\gamma h]
\end{equation}
where we have used the relation $\sqrt{1+2x} \geqs 1+x-\frac12 x^2$. 
\end{proof}


\subsection{Unstable case}
\label{ssec_nu}

We now consider a similar situation as in Section \ref{ssec_ns}, but with
an unstable equilibrium, that is, we assume that $a(t)\geqs a_0>0$ for all
$t\in I$.  Theorem \ref{t_Tihonov} shows the existence of a particular
solution $\xhatdet_t$ of the deterministic equation \eqref{ODEa} such that
$\abs{\xhatdet_t-x^\star(t)} \leqs c_1\eps$ for all $t\in I$. We are
interested in the stochastic process $y_t = x_t - \xhatdet_t$, which
describes the deviation due to noise from this deterministic solution
$\xhatdet$. It obeys the SDE
\begin{equation}
\label{nu1}
\6y_t = \frac1\eps \bigbrak{\ba(t)y_t + \bb(y_t,t)} \6t +
\frac{\sigma}{\sqrt\eps} \6W_t,
\end{equation}
where 
\begin{equation}
\label{nu1b}
\begin{split}
\ba(t) &= \ba_\eps(t) = \sdpar{f}{x}(\xhatdet_t,t) \\
\bb(y,t) &= \bb_\eps(y,t) = f(\xhatdet_t+y,t) - f(\xhatdet_t,t) - \ba(t)y
\end{split}
\end{equation}
are the analogs of $\ba$ and $\bb$ defined in \eqref{ns2}. Taking $\eps$
sufficiently small, we may assume that there exist constants $\ba_0, \ba_1,
\bd > 0$, such that the following estimates hold for all $t\in I$ and all
$y$ such that $\abs{y} \leqs \bd$:
\begin{equation}
\label{nu2}
\ba(t) \leqs -\ba_0, 
\qquad
\abs{\ba'(t)} \leqs \ba_1,
\qquad
\abs{\bb(y,t)} \leqs My^2. 
\end{equation}
The bound on $\abs{\ba'(t)}$ is a consequence of the analog
of~\eqref{ns3d} together with the fact that
$\abs{\xhatdet_0-x^\star(0)}=\Order{\eps}$. 

We first consider the linear equation
\begin{equation}
\label{nu3}
\6y^0_t = \frac1\eps \ba(t) y^0_t \6t + \frac\sigma{\sqrt\eps} \6W_t.
\end{equation}
Given the initial value $y^0_0$, the solution $y^0_t$ at time $t$ is a
Gaussian random variable with mean $y^0_0 \e^{\balpha(t)/\eps}$ and variance
\begin{equation}
\label{nu5}
v(t) = \frac{\sigma^2}\eps \int_0^t \e^{2\balpha(t,s)/\eps} \6s,
\end{equation}
where $\balpha(t,s) = \int_s^t \ba(u)\6u \geqs \ba_0(t-s)$ for $t\geqs s$. 
The variance can be estimated with the help of the following lemma.

\begin{lemma}
\label{l_nu0}
For $0<\eps<2\ba_0^2/\ba_1$, one has
\begin{equation}
\label{nu5a}
\frac1\eps \int_0^{t} \e^{2\balpha(t,s)/\eps} \6s = 
\Bigbrak{\frac{\e^{2\balpha(t)/\eps}}{2\ba(0)} - 
\frac1{2\ba(t)}} 
\bigbrak{1+ \Order{\eps}}.
\end{equation}
\end{lemma}
\begin{proof}
By integration by parts, we obtain that
\begin{equation}
\label{nu5b}
\int_0^t \e^{2\balpha(t,s)/\eps} \6s 
=  \frac{\eps}{2\ba(0)}\e^{2\balpha(t)/\eps} - \frac{\eps}{2\ba(t)} 
- \frac\eps2\int_0^t \frac{\ba'(s)}{\ba(s)^2} 
\e^{2\balpha(t,s)/\eps} \6s,
\end{equation}
which implies that
\begin{equation}
\label{nu5c}
\Bigbrak{1-\frac\eps2\frac{\ba_1}{\ba_0^2}}
\int_0^t \e^{2\balpha(t,s)/\eps} \6s 
\leqs \frac{\eps}{2\ba(0)}\e^{2\balpha(t)/\eps} - \frac{\eps}{2\ba(t)}
\leqs \Bigbrak{1+\frac\eps2\frac{\ba_1}{\ba_0^2}}
\int_0^t \e^{2\balpha(t,s)/\eps} \6s. 
\end{equation}
By our hypothesis on $\eps$, the first term in brackets is positive.
\end{proof}

Unlike in the
stable case, the variance grows exponentially fast (at least with
$\e^{2\ba_0t/\eps}$). If $\rho\geqs\abs{y^0_0}$, we have
\begin{equation}
\label{nu6}
\begin{split}
\bigprobin{0,y^0_0}{\sup_{0\leqs s\leqs t} \abs{y^0_s} < \rho} 
&\leqs \bigprobin{0,y^0_0}{\abs{y^0_t}<\rho} \\
&= \int_{-\rho-y^0_0\e^{\balpha(t)/\eps}}^{\rho-y^0_0\e^{\balpha(t)/\eps}} 
\frac{\e^{-x^2/2v(t)}}{\sqrt{2\pi v(t)}} \6x 
\leqs \frac{2\rho}{\sqrt{2\pi v(t)}},
\end{split}
\end{equation}
which goes to zero as $\rho\sigma^{-1}\e^{-\balpha(t)/\eps}$ for
$t\to\infty$. In this estimate, however, we neglect all trajectories that
leave the interval $(-\rho,\rho)$ before $t$ and come back. We will derive
a more precise estimate for the general, nonlinear case by introducing a
partition of $[0,t]$.

The following proposition, which restates Theorem~\ref{t_unstable} in terms
of $y_t$, is the main result of this subsection.

\begin{prop}
\label{l_nu1}
There exist constants $\eps_0, h_0>0$ such that for all $h\leqs
\sigma\wedge h_0$, all $\eps\leqs\eps_0$ and for any
given $y_0$ with $\abs{y_0}\sqrt{2\ba(0)}< h$, we have 
\begin{equation}
\label{nu7}
\Bigprobin{0,y_0}{\sup_{0\leqs s\leqs t} \abs{y_s} \sqrt{2\ba(s)} < h} 
\leqs \sqrt{\e} 
\exp\Bigset{-\kappa \frac{\sigma^2}{h^2} \frac{\balpha(t)}\eps},
\end{equation} 
where $\kappa = \frac\pi{2{\e}} \bigpar{1-\Order{h}-\Order{\eps}}$.
\end{prop}
\goodbreak
\begin{proof}\hfill
\begin{enum}
\item   Let $K\in\N$ and let $0=u_0 < u_1 < \dots < u_K=t$ be any
partition of the interval $[0,t]$. We define the events 
\begin{equation}
\label{p:lnu:1}
\begin{split}
A_k & = \Bigsetsuch{\w}{\sup_{u_k\leqs s\leqs u_{k+1}} \abs{y_s}
\sqrt{2\ba(s)} < h} \\
B_k & = \Bigsetsuch{\w}{\abs{y_{u_k}}
\sqrt{2\ba(u_k)} < h} \supset A_{k-1}.
\end{split}
\end{equation}
Let $q_k$ be a deterministic upper bound on $P_k =
\probin{u_k,y_{u_k}}{A_k}$, valid on $B_k$. Then we have by the Markov
property
\begin{align}
\label{p:lnu:2}
\nonumber
&\Bigprobin{0,y_0}{\sup_{0\leqs s\leqs t} \abs{y_s} \sqrt{2\ba(s)} < h} \\
\nonumber
&\qquad{} = \Bigprobin{0,y_0}{\bigcap_{k=0}^{K-1} A_k} 
= \Bigexpecin{0,y_0}{\indicator{\bigcap_{k=0}^{K-2} A_k} 
\bigecondin{0,y_0}{\indicator{A_K}}{\set{y_s}_{0\leqs s\leqs u_{K-1}}}} \\
&\qquad{} = \Bigexpecin{0,y_0}{\indicator{\bigcap_{k=0}^{K-2} A_k} P_{K-1}} 
\leqs q_{K-1} \Bigprobin{0,y_0}{\bigcap_{k=0}^{K-2} A_k} 
\leqs \dots \leqs \prod_{k=0}^{K-1} q_k.
\end{align}

\item   To define the partition, we set
\begin{equation}
\label{p:lnu:3}
K = \Bigintpartplus{\frac1\gamma \frac{\balpha(t)}\eps \frac{\sigma^2}{h^2}}
\end{equation}
for some $\gamma\in(0,1]$ to be chosen later, and 
\begin{equation}
\label{p:lnu:4}
\balpha(u_{k+1},u_k) = \gamma\eps\frac{h^2}{\sigma^2},
\qquad k=0,\dots,K-2.
\end{equation}
Since $\balpha(u_{k+1},u_k) \geqs \ba_0(u_{k+1}-u_k)$, we have $u_{k+1}-u_k
\leqs \frac{h^2}{\sigma^2}\frac\gamma{\ba_0}\eps$, and using Taylor's
formula, we find for all $s\in[u_k,u_{k+1}]$ and all $k=0,\dots,K-1$ 
\begin{equation}
\label{p:lnu:5}
1 - \frac{h^2}{\sigma^2}\frac{\ba_1}{\ba_0^2}\gamma\eps 
\leqs \frac{\ba(s)}{\ba(u_k)} \leqs 
1 + \frac{h^2}{\sigma^2}\frac{\ba_1}{\ba_0^2}\gamma\eps,
\end{equation}
where $\ba_1$ is the upper bound on $\abs{\ba'}$, see \eqref{nu2}.
In order to estimate $P_k$, we introduce linear approximations
$(y^{(k)}_t)_{t\in[u_k,u_{k+1}]}$ for $k\in\set{0,\dots,K-2}$, defined
by
\begin{equation}
\label{p:lnu:6}
\6y^{(k)}_t = \frac1\eps \ba(t) y^{(k)}_t + \frac\sigma{\sqrt\eps}
\6W^{(k)}_t, 
\qquad y^{(k)}_{u_k} = y_{u_k},
\end{equation}
where $W^{(k)}_t = W_t - W_{u_k}$ is a Brownian motion with
$W^{(k)}_{u_k}=0$ which is independent of $\setsuch{W_s}{0\leqs s\leqs
u_k}$. If $\w\in A_k$, we have for all $s\in[u_k,u_{k+1}]$
\begin{equation}
\label{p:lnu:7}
\begin{split}
\abs{y_s(\w) - y^{(k)}_s(\w)} &\leqs \frac1\eps \int_{u_k}^s
\e^{\balpha(s,u)/\eps} \abs{\bb(y_u,u)} \6 u \\
&\leqs \frac{Mh^2}{2\ba_0} \frac{\e^{\balpha(u_{k+1},u_k)/\eps}}{\ba(u_k)} 
\bigbrak{1+\Order{\eps}} 
\leqs r_0 \frac{h^2}{\sqrt{2\ba(s)}}, 
\end{split}
\end{equation}
where 
$r_0 = M\e(2 \ba_0^3)^{-1/2} + \Order{\eps}$. This shows that on $A_k$,
\begin{equation}
\label{p:lnu:8}
\abs{y^{(k)}_s(\w)} \leqs \bigbrak{1+r_0 h} \frac h{\sqrt{2\ba(s)}}
\qquad \forall s\in[u_k,u_{k+1}].
\end{equation}

\item   We are now ready to estimate $P_k$. \eqref{p:lnu:8} shows
that on $B_k$,   
\begin{equation}
\label{p:lnu:9}
\begin{split}
P_k &\leqs \Bigprobin{u_k,y_{u_k}}
{\sup_{u_k\leqs s\leqs u_{k+1}} \abs{y^{(k)}_s} \sqrt{2\ba(s)} < h(1+r_0h)}
\\
&\leqs \bigprobin{u_k,y_{u_k}}
{\abs{y^{(k)}_{u_{k+1}}} \sqrt{2\ba(u_{k+1})} < h(1+r_0h)} \\
&\leqs \frac1{\sqrt{2\pi v^{(k)}_{u_{k+1}}}}
\frac{2h(1+r_0h)}{\sqrt{2\ba(u_{k+1})}},
\end{split}
\end{equation}
where $v^{(k)}_{u_{k+1}}$ denotes the conditional variance of
$y^{(k)}_{u_{k+1}}$, given $y_{u_k}$. As in \eqref{nu5a}, 
\begin{equation}
\label{p:lnu:10}
v^{(k)}_{u_{k+1}} = \frac{\sigma^2}\eps \int_{u_k}^{u_{k+1}}
\e^{2\balpha(u_{k+1},s)/\eps} \6s 
= \frac{\sigma^2}2 \biggbrak{\frac{\e^{2\balpha(u_{k+1},u_k)/\eps}}{\ba(u_k)}
- \frac1{\ba(u_{k+1})}} \bigbrak{1+\Order{\eps}}.
\end{equation}
It follows that
\begin{equation}
\label{p:lnu:11}
\begin{split}
\ba(u_{k+1})v^{(k)}_{u_{k+1}} 
&\geqs \frac{\sigma^2}2 \Bigbrak{\e^{2\gamma h^2/\sigma^2}
\frac{\ba(u_{k+1})}{\ba(u_k)}-1} \bigbrak{1-\Order{\eps}} \\
&\geqs \frac{\sigma^2}2 \Bigbrak{\Bigpar{1+2\gamma\frac{h^2}{\sigma^2}}
\Bigpar{1-\frac{\ba_1}{\ba_0^2}\frac{h^2}{\sigma^2}\gamma\eps} - 1}
\bigbrak{1-\Order{\eps}} \\
&\geqs \gamma h^2
\Bigbrak{1-\frac{\ba_1}{2\ba_0^2}\bigpar{1+2\gamma}\eps}
\bigbrak{1-\Order{\eps}} \\ 
&\geqs \gamma h^2 \bigbrak{1-\Order{\eps}}.
\end{split}
\end{equation}
Inserting this into \eqref{p:lnu:9}, we obtain for each $k=0,\dots,K-2$ on
$B_k$ the estimate 
\begin{equation}
\label{p:lnu:12}
P_k \leqs \frac{2h(1+r_0h)}{\sqrt{2\pi}} 
\frac1{\sqrt{2\gamma h^2}}\bigbrak{1+\Order{\eps}} 
= \frac1{\sqrt{\pi\gamma}} 
\bigbrak{1+\Order{\eps}+\Order{h}} \bydef q.
\end{equation}
Note that for any $\gamma\in(1/\pi,1]$, there exist $h_0>0$ and $\eps_0>0$ such
that $q<1$ for all $h\leqs h_0$ and all $\eps\leqs\eps_0$. Since
$q_{K-1}=1$ is an obvious bound, we obtain from \eqref{p:lnu:2} 
\begin{equation}
\label{p:lnu:13}
\Bigprobin{0,y_0}{\sup_{0\leqs s\leqs t} \abs{y_s} \sqrt{2\ba(s)} < h} 
\leqs q^{K-1} 
\leqs \frac1q \exp\Bigset{-\frac{\balpha(t)}{\eps}
\frac{\sigma^2}{h^2} \frac1{2\gamma q^2} q^2\log\bigpar{1/q^2}}.
\end{equation}
Choosing $\gamma$ so that $q^2=1/\e$ holds, yields almost the optimal
exponent, and we obtain 
\begin{equation}
\label{p:lnu:13b}
\Bigprobin{0,y_0}{\sup_{0\leqs s\leqs t} \abs{y_s} \sqrt{2\ba(s)} < h} 
\leqs \sqrt{\e} \exp\Bigset{-\kappa \frac{\balpha(t)}\eps
\frac{\sigma^2}{h^2}}.
\end{equation}
\qed
\end{enum}
\renewcommand{\qed}{}
\end{proof}


\section{Pitchfork bifurcation}
\label{sec_pitchfork}


\subsection{Preliminaries}
\label{ssec_pp}

We consider the nonlinear SDE
\begin{equation}
\label{pp1}
\6 x_\t = \frac{1}{\eps} f(x_\t,\t) \6\t +
\frac{\sigma}{\sqrt{\eps}} \6 W_\t
\end{equation}
in the region $\cM=\setsuch{(x,t)\in\R^2}{\abs{x}\leqs d,\
\abs{t}\leqs T}$. We assume that
\begin{itemiz}
\item   there exists a constant $M>0$ such that $f(x,t)$ is three times
continuously differentiable with respect to $x$ and $t$ and
$\abs{\sdpar{f}{xxx}(x,t)}\leqs 6M$ for all $(x,t)\in\cM$; 
\item   $f(x,t)=-f(-x,t)$ for all $(x,t)\in\cM$;
\item   $f$ exhibits a supercritical pitchfork bifurcation at the
origin, that is (after rescaling),  
\begin{equation} 
\label{pp2} 
\sdpar{f}{x}(0,0) = 0, \qquad
\sdpar{f}{tx}(0,0) = 1  \qquad \text{and} \qquad 
\sdpar{f}{xxx}(0,0) =-6
\end{equation}
\end{itemiz}
Using Taylor series and the symmetry assumptions, we may write for all
$(x,t)\in\cM$ 
\begin{equation}
\label{pp3}
\begin{split}
f(x,t) &= a(t)x + b(x,t) = x\bigbrak{a(t)+g_0(x,t)} \\
\sdpar{f}{x}(x,t) &= a(t) + g_1(x,t)
\end{split}
\end{equation}
where $a(t)$, $g_0(x,t)$, $g_1(x,t)$ are twice continuously differentiable
functions satisfying
\begin{align}
\nonumber
a(t) &= \sdpar{f}{x}(0,t) = t+\Order{t^2} \\
\label{pp4}
g_0(x,t) &= \bigbrak{-1+\gamma_0(x,t)}x^2 & 
\abs{g_0(x,t)} &\leqs  M x^2 \\
\nonumber
g_1(x,t) &= \bigbrak{-3+\gamma_1(x,t)}x^2 & 
\abs{g_1(x,t)} &\leqs 3 M x^2, 
\end{align}
with $\gamma_0, \gamma_1$ some continuous functions such that
$\gamma_0(0,0)=\gamma_1(0,0)=0$.
The following standard result from bifurcation theory is easily obtained by
applying the implicit function theorem, see~\cite[p.~150]{GH}
or~\cite[Section~II.4]{IJ} for instance. We state it without proof. 

\begin{prop}
\label{p_pp1}
If $T$ and $d$ are sufficiently small, there exist twice continuously
differentiable functions $x^\star,\bar x: (0,T]\to\R_+$ of the form
\begin{equation}
\label{pp5}
\begin{split}
x^\star(t) &= \sqrt t \bigbrak{1+\orderone{T}} \\
\bar x(t) &= \sqrt{t/3} \bigbrak{1+\orderone{T}}
\end{split}
\end{equation}
with the following properties:
\begin{itemiz}
\item   the only solutions of $f(x,t)=0$ in $\cM$ are either of the form
$(0,t)$, or of the form $(\pm x^\star(t),t)$ with $t>0$;
\item   the only solutions of $\sdpar fx(x,t)=0$ in $\cM$ are of the form
$(\pm\bar x(t),t)$ with $t\geqs 0$;
\item   the derivative of $f$ at $\pm x^\star(t)$ is 
\begin{equation}
\label{pp6}
a^\star(t) = \sdpar fx(x^\star(t),t) = -2t\bigbrak{1+\orderone{T}}.
\end{equation}
\item   the derivatives of $x^\star(t)$ and $\bar x(t)$ satisfy
\begin{equation}
\label{pp6b}
\dtot{x^\star}t = \frac1{2\sqrt t}[1+\orderone{T}], \qquad
\dtot{\bar x}t = \frac1{2\sqrt{3t}}[1+\orderone{T}].
\end{equation}
\end{itemiz}
\end{prop}

As already pointed out in Section~\ref{ssec_rb}, there is no
restriction in assuming $T$ and $d$ to be small. Thus we may assume
that the terms $\orderone{T}$ are sufficiently small to do no
harm. For instance, we may and will always assume that $a^\star(t)<0$.

Equation \eqref{pp4} also implies the existence of constants $a_+\geqs
a_->0$ such that 
\begin{equation}
\label{pp7}
\begin{split}
a_+ t \leqs a(t) \leqs a_- t \qquad &\text{for $-T\leqs t\leqs 0$} \\
a_- t \leqs a(t) \leqs a_+ t \qquad &\text{for $0\leqs t\leqs T$.}
\end{split}
\end{equation} 
The function $\alpha(t,s)=\int_s^t a(u)\6 u$ thus satisfies
\begin{align}
\nonumber
-\tfrac12 a_+(s^2-t^2) &\leqs \alpha(t,s) \leqs -\tfrac12 a_-(s^2-t^2) 
& &\text{if $s\leqs t\leqs 0$} \\
\label{pp8}
\tfrac12 a_-t^2 - \tfrac12 a_+s^2 &\leqs \alpha(t,s) 
\leqs \tfrac12 a_+ t^2 - \tfrac12 a_- s^2
& &\text{if $s\leqs 0\leqs t$} \\
\nonumber
\tfrac12 a_-(t^2-s^2) &\leqs \alpha(t,s) \leqs \tfrac12 a_+(t^2-s^2) 
& &\text{if $0\leqs s\leqs t$.} 
\end{align}

We are going to analyse the dynamics in three different regions of the
$(t,x)$-plane: near $x=0$ for $t\leqs\sqrt\eps$, near $x=0$ for
$t\geqs\sqrt\eps$, and near $x=x^\star(t)$ for $t\geqs\sqrt\eps$. In order
to delimit the last two regions, we introduce (somewhat arbitrarily) the
function  
\begin{equation}
\label{pp13}
\tx(t) = \sqrt\lambda \, x^\star(t),
\end{equation}
set 
\begin{equation}
\label{pp14}
\ta(t) = \sdpar fx(\tx(t),t),
\end{equation}
and define the region
\begin{equation}
\label{pp15}
\cD = \bigsetsuch{(x,t)}
{\sqrt\eps \leqs t \leqs T, \abs{x} < \tx(t)}, 
\end{equation}
which has the following properties:
\begin{enum}
\item[(a)]      for all $(x,t)\in\cD$ with $x\ne0$, one has
\begin{equation}
\label{pp16a}
\frac1x f(x,t) \geqs \kappa a(t) \qquad
\text{with $\kappa = 1 - \lambda - \orderone{T}$.} 
\end{equation}
\item[(b)]  for all $(x,t)\in[-d,d\mskip2mu]\times[\sqrt\eps,T]\setminus\cD$, 
\begin{equation}
\label{pp16b}
\sdpar fx(x,t) \leqs \ta(t) \leqs -\mm a(t) \qquad
\text{with $\mm = 3\lambda - 1 - \orderone{T}$.}
\end{equation}
\end{enum}
For our results to be of interest, $\kappa>0$ and $\mm>0$ are
necessary, which requires $\lambda\in(\tfrac13,1)$. As we shall see,
we will actually need $\lambda\in (\tfrac13,\tfrac12)$. Furthermore, in
Section~\ref{ssec_pesc}, we need to assume that
$\sigma\abs{\log\sigma}^{3/2}=\Order{\sqrt\eps}$.

In the following subsections, we investigate the three different regimes: 
In Section~\ref{ssec_psmallt}, we analyse the behaviour for
$t\leqs\sqrt\eps$. Theorem~\ref{t_before} is proved in the same way as
Theorem~\ref{t_stable}, the main difference lying in the behaviour of the
variance which is investigated in Lemma~\ref{l_ps1}.
 
Section~\ref{ssec_pesc} is devoted to the rather involved proof of 
Theorem~\ref{t_escape}. We start by giving some preparatory
results. Proposition~\ref{l_escape} estimates the probability of
remaining in a smaller strip $\cS$ in a similar way as
Proposition~\ref{l_nu1}. We then show in Lemma~\ref{l_ez} that the
paths are likely to leave $\cD$ as well, unless the solution of a
suitably chosen linear SDE returns to zero. The probability of such a
return to zero is studied in Lemma~\ref{l_rz}. Finally,
Theorem~\ref{t_escape} is proved, the proof being based on an
iterative scheme.

The last subsection analyses the motion after $\tau_\cD$. Here, the main
difficulty is to control the behaviour of the deterministic solutions,
which are shown to approach $x^\star(t)$, cf.\
Proposition~\ref{p_app1}. We then prove that the paths of the random
process are likely to stay in a neighbourhood of the deterministic
solutions. The proof is similar to the corresponding proof in
Section~\ref{ssec_ns}. 


\subsection{The behaviour for $t\leqs\sqrt\eps$}
\label{ssec_psmallt}

We first consider the linear equation
\begin{equation}
\label{ps1}
\6 x^0_t = \frac1\eps a(t) x^0_t \6 t + \frac\sigma{\sqrt\eps} \6 W_t
\end{equation}
with initial condition $x^0_{t_0} = x_0$ at time $t_0\in[-T,0)$. Let
\begin{equation}
\label{ps3}
v(t,t_0) = \frac{\sigma^2}\eps \int_{t_0}^t \e^{2\alpha(t,s)/\eps} \6 s.
\end{equation}
denote the variance of $x^0_t$. As before, we now introduce a function
$\z(t)$ which will allow us to define a strip that the process $x_t$
is unlikely to leave before time $\sqrt\eps$, see
Corollary~\ref{c_ps1} below. Let 
\begin{equation}
\label{ps4}
\z(t) = \frac1{2\abs{a(t_0)}} \e^{2\alpha(t,t_0)/\eps} + 
\frac1\eps\int_{t_0}^t \e^{2\alpha(t,s)/\eps} \6 s.
\end{equation}
The following lemma describes the behaviour of $\z(t)$.

\begin{lemma}
\label{l_ps1}
Assuming $\eps\leqs 4a(t_0)^2\wedge(t_0/2)^2$, there exist constants $c_\pm
= c_\pm(a_+,a_-)$ such that
\begin{align}
\nonumber
\frac{c_-}{\abs t} &\leqs \z(t) \leqs \frac{c_+}{\abs{t}} 
&&\text{for $t_0\leqs t\leqs -\sqrt\eps$} \\
\label{ps5}
\frac{c_-}{\sqrt\eps} &\leqs \z(t) \leqs \frac{c_+}{\sqrt\eps} 
&&\text{for $-\sqrt\eps\leqs t\leqs \sqrt\eps$} \\
\nonumber
\frac{c_-}{\sqrt\eps}\e^{2\alpha(t)/\eps} &\leqs \z(t) \leqs
\frac{c_+}{\sqrt\eps}\e^{2\alpha(t)/\eps}  
&&\text{for $\sqrt\eps\leqs t\leqs T$.} 
\end{align}
If, moreover, $a'(t)>0$ on $[t_0,t]$, then $\z(t)$ is increasing on
$[t_0,t]$.
\end{lemma}

\begin{proof}
The upper bounds are easy to obtain.
For $t_0\leqs t\leqs -\sqrt\eps$ we have, using $t^2-s^2\leqs 2t(t-s)$,
\begin{equation}
\label{p:lps1}
\z(t)\leqs \frac1\eps \int_{t_0}^t \e^{a_-(t^2-s^2)/\eps}\6 s 
+ \frac1{2\abs{a(t_0)}} \leqs \frac1{\abs t}
\Bigbrak{\frac1{2a_-}+\frac1{2a_+}}.
\end{equation}
For $-\sqrt\eps\leqs t\leqs 0$, the hypothesis $\eps\leqs 4a(t_0)^2$ implies
\begin{equation}
\label{p:lps2}
\z(t)\leqs \frac1\eps \e^{-a_-}\int_{t_0}^0 \e^{-a_-s^2/\eps}\6 s 
+ \frac1{2\abs{a(t_0)}} \leqs \frac1{\sqrt\eps}
\Bigbrak{\e^{-a_-}\int_{-\infty}^0 \e^{-a_-u^2}\6 u + 1}.
\end{equation}
For $0\leqs t\leqs\sqrt\eps$, a similar estimate is obtained by
splitting the integrals for $s\leqs0$ and $s\geqs0$. 
For $t\geqs\sqrt\eps$, we have
\begin{equation}
\label{p:lps3}
\e^{-2\alpha(t)/\eps}\z(t) \leqs \frac1{\sqrt\eps}
\Bigbrak{\int_{-\infty}^0 \e^{-a_-u^2} \6 u + \int_0^{\infty}
\e^{-a_+u^2}\6 u + 1}.
\end{equation}
To obtain the lower bound, we first consider the interval $t_0\leqs t\leqs
\frac12 t_0$, where we use the estimate $t^2-s^2\geqs 2t_0(t-s)$,
valid for all $s\in[t_0,t]$, which yields
\begin{equation}
\label{p:lps4}
\z(t)\geqs \frac1\eps \int_{t_0}^t\e^{-2a_+\abs{t_0}(t-s)/\eps}\6 s +
\frac{\e^{-2a_+\abs{t_0}(t-t_0)/\eps}}{2a_+\abs{t_0}} \geqs 
\frac1{2a_+\abs t}. 
\end{equation}
For $\frac12 t_0\leqs t\leqs -\sqrt\eps$, we have $t^2-s^2\geqs 3t(t-s)$
for all $s\in[2t,t]$, and thus
\begin{equation}
\label{p:lps5}
\z(t)\geqs \frac1\eps \int_{2t}^t\e^{-3a_+\abs{t}(t-s)/\eps}\6 s 
\geqs \frac{1-\e^{-3a_+}}{3a_+\abs{t}},
\end{equation}
where we used the relation $t_0\leqs-2\sqrt\eps$ in the last step. By
the same relation, we obtain 
\begin{align}
\label{p:lps6}
\z(t) &\geqs \frac1{\sqrt\eps}\int_{-2}^{-1} \e^{-a_+u^2}\6 u 
&&\text{for $-\sqrt\eps\leqs t\leqs\sqrt\eps$,} \\
\label{p:lps7}
\e^{-2\alpha(t)/\eps}\z(t) &\geqs \frac1{\sqrt\eps}\int_0^1 \e^{-a_+u^2}\6 u 
&&\text{for $t\geqs\sqrt\eps$.} 
\end{align}
Finally, assume that $a'(t)>0$ for all $t$, and recall that $\z(t)$ is
the solution of the initial value problem
\begin{equation}
\label{p:lps8}
\dtot{\z}{t} = \frac{2a(t)}{\eps} \z + \frac1\eps, 
\qquad \z(t_0) = \frac1{2\abs{a(t_0)}}.
\end{equation}
Since $\z(t)\geqs 0$, $\z'>0$ for all positive $t$. For negative
$t$, $\z'$ is positive whenever the function $V(t) = \z(t) +
1/2a(t)$ is negative. We have $V(t_0)=0$ and  
\begin{equation}
\label{p:lps9}
\dtot{V}{t} = \frac{2a(t)}{\eps} V - \frac{a'(t)}{2a(t)^2}.
\end{equation}
Since $V'<0$ whenever $V=0$, $V$ can never become positive. This
implies $\z'\geqs 0$. 
\end{proof}

The following proposition shows that the solution $x^0_t$ of the linearized
equation~\eqref{ps1} is likely to track the solution of the
corresponding deterministic equation.

\begin{prop}
\label{l_ps2}
Assume that $-T\leqs t_0<t\leqs\sqrt\eps$. For sufficiently small $\eps$, 
\begin{equation}
\label{ps6}
\Bigprobin{t_0,x_0}{\sup_{t_0\leqs s\leqs t}
\frac{\abs{x^0_s-x_0\e^{\alpha(s,t_0)/\eps}}}{\sqrt{\z(s)}} > h}
\leqs C(t,\eps)
\exp\Bigset{-\frac12\frac{h^2}{\sigma^2}\bigbrak{1-r(\eps)}},
\end{equation}
where 
\begin{equation}
\label{ps7}
C(t,\eps) = \frac{\abs{\alpha(t,t_0)}}{\eps^2} + \frac{a_++4\sqrt\eps+4}\eps
\end{equation}
and where $r(\eps)=\Order{\eps}$ for $t_0\leqs t\leqs -\sqrt\eps$,
and $r(\eps)=\Order{\sqrt\eps}$ for $-\sqrt\eps\leqs t\leqs\sqrt\eps$. 
\end{prop}

\begin{proof}
Let $t_0 = u_0 < \dots < u_K=t$ be a partition of the interval
$[t_0,t]$. By Lemma~\ref{l_ns2}, the probability in~\eqref{ps6} is
bounded by $2\sum_{k=1}^K P_k$, where 
\begin{equation}
\label{p:lps2:1}
P_k = \exp\Bigset{-\frac12\frac{h^2}{\sigma^2} \frac1{\z(u_k)}
\inf_{u_{k-1}\leqs u\leqs u_k} \z(u) \e^{2\alpha(u_k,u)/\eps}}.
\end{equation}
If $t\leqs-\sqrt\eps$, we define the partition by
\begin{equation}
\label{p:lps2:2}
K = \biggintpartplus{\frac{-\alpha(t,t_0)}{2\eps^2}}, \qquad
-\alpha(u_k,t_0) = 2\eps^2 k \quad
\text{for $k=0,\dots,K-1$}.
\end{equation}
Estimating $P_k$ as in the proof of Proposition~\ref{l_ns3}, we obtain
\begin{equation}
\label{p:lps2:5}
P_k \leqs \exp\Bigset{-\frac12\frac{h^2}{\sigma^2}
\Bigpar{1-\frac{2\eps}{a_-c_-}}\e^{-4\eps}}.
\end{equation}
Therefore, \eqref{ps6} holds with $C(t,\eps)=\abs{\alpha(t,t_0)}/\eps^2+2$.

For $-\sqrt\eps\leqs t\leqs\sqrt\eps$, we define the partition
separately in two different regions. Let
\begin{equation}
\label{p:lps2:6}
K_0 = \biggintpartplus{\frac{-\alpha(-\sqrt\eps,t_0)}{2\eps^2}}, \qquad
K = K_0 + \biggintpartplus{\frac{t+\sqrt\eps}\eps}.
\end{equation}
The partition times are defined via
\begin{align}
\nonumber
-\alpha(u_k,t_0)  &= 2\eps^2 k && 
\text{for $0\leqs k\leqs K_0-1$} \\
\label{p:lps2:7}
u_k &= -\sqrt\eps + \eps(k-K_0) && 
\text{for $K_0\leqs k\leqs K-1$}.
\end{align}
In the first case, we immediately obtain the bound \eqref{p:lps2:5}.
In the second case, estimating $P_k$ in the usual way shows that 
\begin{equation}
\label{p:lps2:10}
P_k \leqs \exp\Bigset{-\frac12\frac{h^2}{\sigma^2}
\Bigpar{1-\frac{\sqrt\eps}{c_-}[1+2a_+c_+]}\e^{-a_+\eps}}.
\end{equation}
Finally, let us note that, for $-\sqrt\eps\leqs t\leqs \sqrt\eps$,
\begin{equation}
\label{p:lps2:11new}
2K 
\leqs \frac{\abs{\alpha(-\sqrt\eps,t_0)}}{\eps^2} +
\frac2\eps(t+\sqrt\eps) + 4
\leqs \frac{\abs{\alpha(t,t_0)}}{\eps^2} +\frac{a_+}\eps
+\frac4{\sqrt\eps} + 4,
\end{equation}
which concludes the proof of the proposition.
\end{proof}

Let us now compare solutions of the two SDEs
\begin{align}
\label{ps8a}
\6 x^0_t &= \frac1\eps a(t) x^0_t \6 t + \frac\sigma{\sqrt\eps} \6 W_t &
x^0_{t_0} &= x_0 \\
\label{ps8b}
\6 x_t &= \frac1\eps f(x_t,t) \6 t + \frac\sigma{\sqrt\eps} \6 W_t &
x_{t_0} &= x_0, 
\end{align}
where $t_0\in[-T,0)$.
We define the events
\begin{align}
\label{ps9a}
\Omega^0_t(h) &= \Bigsetsuch{\w}{\bigabs{x^0_s(\w)-x_0\e^{\alpha(s,t_0)/\eps}} 
\leqs h\sqrt{\z(s)} \; \forall s \in [t_0,t]} \\
\label{ps9b}
\Omega_t(h) &= \Bigsetsuch{\w}{\bigabs{x_s(\w)-x_0\e^{\alpha(s,t_0)/\eps}} 
\leqs h\sqrt{\z(s)} \;\forall s \in [t_0,t]}. 
\end{align}
Proposition~\ref{l_ps2} gives us an upper bound on the probability of the
complement of $\Omega^0_t(h)$. We now give relations between these events.

\begin{prop}
\label{l_ps3}
Let $t\in[t_0,\sqrt\eps\mskip3mu]$ and $\abs{x_0}\leqs
h/\eps^{1/4}$, where we assume $h^2<\eps/\gamma$ 
for $\gamma=M(1+2\sqrt{c_+})^3c_+/\sqrt{c_-}$ and $h^2\leqs
d^2\sqrt\eps/(1+2\sqrt{c_+})^2$. Then
\begin{align}
\label{ps10a}
\Omega_t(h) 
&\subas \Omega^0_t\Bigpar{\Bigbrak{1+\gamma\frac{h^2}\eps} h} \\
\label{ps10b}
\Omega^0_t(h) 
&\subas \Omega_t\Bigpar{\Bigbrak{1+\gamma\frac{h^2}\eps } h}. 
\end{align}
\end{prop}

\begin{proof}
Assume first that $\w\in\Omega^0_t(h)$ and let $\delta = \gamma
h^2/\eps$. Then we have $\delta<1$ by assumption. By~\eqref{pp3}, the
difference $z_s=x_s-x^0_s$ satisfies 
\begin{equation}
\label{p:lps3:1}
z_s = \frac1\eps \int_{t_0}^s \e^{\alpha(s,u)/\eps} b(x_u,u) \6 u.
\end{equation}
We consider the first exit time 
\begin{equation}
\label{p:lps3:3}
\tau = \inf \bigsetsuch{s\in[t_0,t]}{\abs{z_s} \geqs \delta h\sqrt{\z(s)}} 
\in [t_0,t] \cup \set\infty.
\end{equation}
For all $\w$ in the set 
\begin{equation}
\label{p:lps3:4}
A = \Omega^0_t(h) \cap \bigsetsuch{\w}{\tau(\w) < \infty},
\end{equation}
and $s\in[t_0,\tau(\w)]$, we have by the hypotheses on $h$ and $x_0$
together with Lemma \ref{l_ps1}
\begin{equation}
\label{p:lps3:neu} 
\abs{x_s(\w)}\leqs \abs{x_0}+ h\sqrt{\z(s)}
\leqs \bigpar{1+(1+\delta)\sqrt{c_+}}\frac h{\eps^{1/4}}\leqs d.
\end{equation}
Therefore, \eqref{pp4} yields
\begin{equation}
\label{p:lps3:2neu}
\abs{z_s} \leqs M\Bigbrak{\bigpar{1+(1+\delta)\sqrt{c_+}}
\frac h{\eps^{1/4}}}^3\;
\frac1\eps \int_{t_0}^s \e^{\alpha(s,u)/\eps} \6 u.
\end{equation}
The integral is bounded by $2\z_{2\eps}(s)$, which can be estimated by
Lemma \ref{l_ps1} once again. Thereby, we obtain
\begin{equation}
\label{p:lps3:5}
\abs{z_s} \leqs M\bigpar{1+(1+\delta)\sqrt{c_+}}^3
\frac{c_+}{\sqrt{c_-}} \frac{h^2}\eps h \sqrt{\z(s)} 
< \delta h \sqrt{\z(s)}, 
\end{equation}
which leads to a contradiction for $s=\tau(\w)$. We conclude that
$\fP(A) = 0$, and thus $\tau(\w)=\infty$ for $\fP$-almost all
$\w\in\Omega^0_t(h)$. This shows that $\abs{z_s(\w)}<\delta h\sqrt{\z(s)}$
and thus $\abs{x_s(\w)-x_0\e^{\alpha(s,t_0)/\eps}}<(1+\delta)
h\sqrt{\z(s)}$ for all these $\w$ and all $s\in[t_0,t]$, which proves
\eqref{ps10b}. 
The proof of the inclusion \eqref{ps10a} is straightforward, using the
same estimates.
\end{proof}

The two preceding propositions immediately imply the main result on the
behaviour of the solution of the nonlinear equation~\eqref{ps8b} for
$t\leqs\sqrt\eps$, i.e., Theorem~\ref{t_before}, which we restate here
with an arbitrary initial time $t_0\in[-T,\sqrt\eps\mskip3mu]$.

\begin{cor}
\label{c_ps1}
Assume that $-T\leqs t_0<t\leqs\sqrt\eps$. Then there exists an
$h_0>0$ such that for all $h\leqs h_0\sqrt\eps$ and all initial
conditions $x_0$ with $\abs{x_0}\leqs h/\eps^{1/4}$, the following
estimate holds:
\begin{equation}
\label{ps6nonlin}
\Bigprobin{t_0,x_0}{\sup_{t_0\leqs s\leqs t}
\frac{\abs{x_s-x_0\e^{\alpha(s,t_0)/\eps}}}{\sqrt{\z(s)}} > h}
\leqs C(t,\eps)
\exp\Bigset{-\frac12\frac{h^2}{\sigma^2}\bigbrak{1-r(\eps)-\Order{h^2/\eps}}},
\end{equation}
where $C(t,\eps)$ and $r(\eps)$ are given in Proposition~\ref{l_ps2}.
\end{cor}


\subsection{Escape from the origin}
\label{ssec_pesc}

We now consider the SDE \eqref{pp1}, written in the form
\begin{equation}
\label{pe1}
\6x_t = \frac1\eps \bigbrak{a(t)x_t + b(x_t,t)} \6t + 
\frac\sigma{\sqrt\eps} \6W_t,
\end{equation}
for $t\geqs t_0\geqs\sqrt\eps$, where we assume that
$\abs{x_{t_0}}\leqs\tx(t_0)$. Our aim is to estimate the
first exit time $\tau_\cD$ of $x_t$ from $\cD$ defined in
\eqref{pp15}. We recall that $a(t) + \frac1x b(x,t) \geqs \kappa a(t)$
in $\cD$, see \eqref{pp16a}. Moreover, we have $a_-t\leqs a(t)\leqs
a_+t$, $0\leqs a'(t)\leqs a_1$, and $\abs{b(x,t)}\leqs M\abs{x}^3$ in
$\cD$.  

We first state a result allowing to estimate the variance of the
linearization of \eqref{pe1}.

\begin{lemma}
\label{l_unst_var}
Let $a(t)$ be any continuously differentiable, strictly positive, increasing
function, and set $\alpha(t,s)=\int_s^t a(u)\6 u$.
Then the integral
\begin{equation}
\label{luv1}
v(t,s) = \frac{\sigma^2}{\eps} \int_s^t \e^{2\alpha(t,u)/\eps} \6 u
\end{equation}
satisfies the inequalities
\begin{equation}
\label{luv2}
\frac{\sigma^2}{2a(t)} \bigbrak{\e^{2\alpha(t,s)/\eps}-1} 
\leqs  v(t,s) \leqs \frac{\sigma^2}{2a(s)}\e^{2\alpha(t,s)/\eps}. 
\end{equation}
\end{lemma}
\begin{proof}
Using integration by parts, we have 
\begin{equation}
\label{pluv:1}
\e^{-2\alpha(t,s)/\eps} v(t,s) = 
\sigma^2 \Bigbrak{\frac1{2a(s)} - \frac1{2a(t)}\e^{-2\alpha(t,s)/\eps}
- \int_s^t \frac{a'(u)}{2a(u)^2} \e^{-2\alpha(u,s)/\eps} \6 u}.
\end{equation}
The upper bound follows immediately, and the lower bound is obtained by 
bounding the exponential in the last integral by $1$.
\end{proof}

Our first step towards estimating $\tau_\cD$ is to estimate the first exit
time $\tau_\cS$ from a smaller strip $\cS$, defined as
\begin{equation}
\label{pe2}
\cS = \biggsetsuch{(x,t)}
{\sqrt\eps\leqs t\leqs T, \abs{x} < \frac h{\sqrt{a(s)}}},
\end{equation}
where we will choose 
\begin{equation}
\label{pe2b}
h = 2 \sigma \sqrt{\abs{\log\sigma}}.
\end{equation}

\begin{prop}
\label{l_escape}
Let $t_0\geqs\sqrt\eps$ and $\abs{x_0}\leqs h/\sqrt{a(t_0)}$. Then,
for any $\mu>0$, we have 
\begin{equation}
\label{esc2}
\bigprobin{t_0,x_0}{\tau_\cS \geqs t} 
\leqs \Bigpar{\frac h\sigma}^\mu
\exp\biggset{-\frac{\mu}{1+\mu}\frac{\alpha(t,t_0)}{\eps} 
\Bigbrak{1-\BigOrder{\frac1{\mu\log(h/\sigma)}}}}
\end{equation}
under the condition
\begin{equation}
\label{esc3}
\Bigpar{\frac h\sigma}^{3+\mu} 
\BigOrder{\log\frac{h}{\sigma}} \leqs \frac{t_0^2}{\sigma^2}.
\end{equation}
\end{prop}
\goodbreak
\begin{proof}\hfill
\begin{enum}
\item   For $K\in\N$, we introduce a partition $t_0=u_0<\dots<u_K=t$
of the interval $[t_0,t]$, which will be chosen later, and for each
$k$, we define a linear approximation
$(x_t^{(k)})_{t\in[u_k,u_{k+1}]}$ by  
\begin{equation}
\label{pesc:2}
\6 x^{(k)}_t = \frac1\eps a(t)x^{(k)}_t \6 t 
+ \frac{\sigma}{\sqrt\eps} \6 W^{(k)}_t
\qquad x^{(k)}_{u_k}=x_{u_k},
\end{equation}
where $W^{(k)}_t = W_t - W_{u_k}$.
Assume that $\abs{x_s}\sqrt{a(s)}\leqs h$ for all $s\in[u_k,u_{k+1}]$. Then 
by Lemma \ref{l_unst_var} 
\begin{equation}
\label{pesc:3}
\begin{split}
\abs{x_s-x^{(k)}_s} 
&\leqs \frac1\eps \int_{u_k}^s \abs{b(x_u,u)} \e^{\alpha(s,u)/\eps} \6 u \\
&\leqs M \frac{h^3}{a(u_k)^{3/2}} \frac1{a(u_k)}
\e^{\alpha(u_{k+1},u_k)/\eps} \leqs \frac{h}{\sqrt{a(s)}} 
\end{split}
\end{equation}
for  $s\in[u_k,u_{k+1}]$, provided the partition is chosen in such a
way that for all $k$ 
\begin{equation}
\label{pesc:1}
h^2 \leqs \frac{a_-^2}M \sqrt{\frac{a(u_k)}{a(u_{k+1})}}
\e^{-\alpha(u_{k+1},u_k)/\eps} t_0^2.
\end{equation}

\item   If $\abs{x_{u_k}} \sqrt{a(u_k)} \leqs h$, then we have 
\begin{align}
\nonumber
\Bigprobin{u_k,x_{u_k}}{\sup_{u_k\leqs s\leqs u_{k+1}}
\abs{x_s}\sqrt{a(s)} \leqs h}
&\leqs \Bigprobin{u_k,x_{u_k}}{\abs{x^{(k)}_{u_{k+1}}}\sqrt{a(u_{k+1})} \leqs
2h} \\
\label{pesc:4}
&\leqs \frac{4h}{\sqrt{2\pi v^{(k)}_{u_{k+1}} a(u_{k+1})}},
\end{align}
where the variance
\begin{equation}
\label{pesc:5}
v^{(k)}_{u_{k+1}} = \frac{\sigma^2}{\eps} \int_{u_k}^{u_{k+1}}
\e^{2\alpha(u_{k+1},s)/\eps} \6 s
\end{equation}
can be estimated by Lemma \ref{l_unst_var}.
We thus have by the Markov property
\begin{equation}
\label{pesc:6}
P = \Bigprobin{t_0,x_0}{\sup_{t_0\leqs s\leqs t}
\abs{x_s}\sqrt{a(s)} \leqs h} 
\leqs \prod_{k=0}^{K-1} \biggpar{\frac4{\sqrt{2\pi}}
\frac{h}{\sqrt{v^{(k)}_{u_{k+1}} a(u_{k+1})}}\wedge 1}.
\end{equation}

\item   We now choose the $u_k$ in such a way that
$v^{(k)}_{u_{k+1}}a(u_{k+1})$ is approximately constant. Given $\mu>0$,
let
\begin{equation}
\label{pesc:7}
\ell = \frac8\pi h^2 \Bigpar{\frac{h^2}{\sigma^2}}^\mu
\end{equation}
(Observe that $\ell \geqs 8h^2/\pi > \sigma^2/2$.) 
Choosing $K$ as the smallest integer satisfying
\begin{equation}
\label{pesc:9}
K \geqs \frac{2\alpha(t,t_0)}{\eps\log(2\ell/\sigma^2)},
\end{equation}
we define the partition by the relations
\begin{align}
\label{pesc:8}
\alpha(u_{k+1},u_k) &{}= \frac\eps2 \log \frac{2\ell}{\sigma^2},
\qquad \text{for $k\in\{0,\dots,K-2\}$,}\\
0<\alpha(u_K,u_{K-1})&{}\leqs \frac\eps2 \log \frac{2\ell}{\sigma^2}.
\end{align}
Then we have 
\begin{equation}
\label{pesc:10}
P \leqs \biggpar{\frac4{\sqrt\pi} \frac h\sigma
\frac1{\sqrt{2\ell/\sigma^2-1}}}^{K-1} 
\leqs \Bigpar{\frac h\sigma}^\mu
\exp\biggset{-\frac{\alpha(t,t_0)}\eps
\frac{\log\bigbrak{\bigpar{\frac{h^2}{\sigma^2}}^\mu - \frac\pi{16}
\frac{\sigma^2}{h^2}}}{\log\bigbrak{\frac{16}\pi
\bigpar{\frac{h^2}{\sigma^2}}^{1+\mu}-1}}},
\end{equation}
which proves \eqref{esc2}. 

\item   It remains to show that condition \eqref{pesc:1} is satisfied. 
Since
\begin{equation}
\label{pesc:11}
\frac{a(u_{k+1})}{a(u_k)} \leqs 1+\frac{a_1}{a(u_k)} (u_{k+1}-u_k) 
\leqs 1 + \frac{a_1\eps}{2a_-^2 t_0^2} \log\biggl\{\frac{16}\pi
\Bigpar{\frac{h^2}{\sigma^2}}^{1+\mu}\biggr\},
\end{equation}
the condition reduces to
\begin{equation}
\label{pesc:12}
\Bigpar{\frac h\sigma}^{3+\mu} 
\biggpar{1 + \frac{a_1}{4a_-^2}\frac{\eps}{t_0^2} \log\Bigbrak{\frac{16}\pi
\Bigpar{\frac{h^2}{\sigma^2}}^{1+\mu}}} \leqs \frac{a_-^2}M \frac{\sqrt\pi}4
\frac{t_0^2}{\sigma^2},
\end{equation}
which is satisfied whenever condition \eqref{esc3} is satisfied.\qed
\end{enum}
\renewcommand{\qed}{}
\end{proof}

We want to choose $\mu$ in such a way that 
$\probin{t_0,x_0}{\tau_\cS \geqs t} \leqs (h/\sigma)^\mu
\e^{-\kappa\alpha(t,t_0)/\eps}$ holds with the same $\kappa$ as in
\eqref{pp16a}. We opt for $\mu=2$, because this choice guarantees the
above estimate for all possible $\kappa$ without choosing a
$\kappa$-dependent $\mu$. For $h=2\sigma\sqrt{\abs{\log\sigma}}$,
Condition~\eqref{esc3} becomes a consequence of the following slightly
stronger condition 
\begin{equation}
\label{pe3}
\sigma \abs{\log\sigma}^{3/2} = \Order{\sqrt\eps},
\end{equation}
which we will assume to be satisfied from now on for the rest of this
subsection. 

\begin{figure}
 \centerline{\psfig{figure=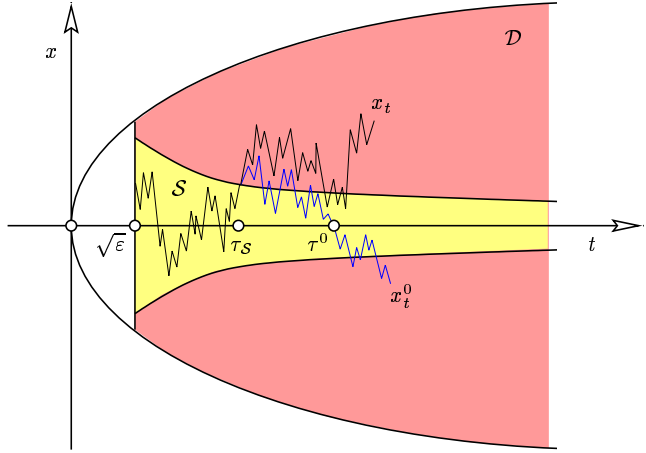,height=60mm,clip=t}}
 \caption[]
 {Assume the path $x_t$ exits the region $\cS$ at time $\tau_\cS$, say
 by passing through the upper boundary of $\cS$. We introduce a
 process $x^0_t$, starting on the same boundary at time $\tau_\cS$,
 which obeys the linear SDE~\eqref{pe4}. Let $\tau^0$ be the time of
 first return to zero of $x^0_t$. Then $x_t$ lies above $x^0_t$ for
 $\tau_\cS<t\leqs\tau^0$. In case $x_t$ also becomes negative,
 the two processes may cross each other. The probability of $x^0_t$ ever
 returning to zero is bounded by $\sigma^{4\kappa}$. If $x^0_t$ does not
 return to zero, $x_t$ is likely to leave $\cD$.}
\label{fig3}
\end{figure}

The second step is to control the probability that $x_t$ returns to zero
after it has left the strip $\cS$. To do so, we will compare solutions of
\eqref{pe1} with those of the linear equation
\begin{equation}
\label{pe4}
\6x^0_t = \frac1\eps a_0(t) x^0_t \6t + \frac\sigma{\sqrt\eps}\6W_t,
\end{equation}
where $a_0(t) = \kappa a(t)$ satisfies $ a_0(t)\leqs f(x,t)/x$ in $\cD$.
The following lemma shows that this choice of $a_0(s)$ implies that
$\abs{x_s}\geqs\abs{x^0_s}$ holds as long as $x_s$ does not return to zero 
(\figref{fig3}).
This implies that if $x^0_s$ does not return to zero before time $t$, then
$x_s$ is likely to leave $\cD$ before time $t$ without returning to zero.

\begin{lemma}
\label{l_ez}
Let $t_0\geqs\sqrt\eps$ and assume that $0<x_0<\tx(t_0)$. We define
\begin{equation}
\label{ez1}
{\cD}^+(t) = \bigsetsuch{(x,s)}
{\text{$\sqrt\eps\leqs s\leqs t$ and $0<x<\tx(s)$}}
\end{equation}
and denote by $\tau_{{\cD}^+}$ the first exit time of $x_s$ from
${\cD}^+(t)$.
Let $\tau^0$ be the time of first return to zero of $x^0_s$ in $[t_0,t]$,
where we set $\tau^0=\infty$ if $x^0_s>0$ for all $t\in[t_0,t]$. Then
$x_s\geqs x^0_s$ for all $s\leqs \tau_{\cD^+}\wedge t$ and
\begin{equation}
\label{ez3}
\begin{split}
\Bigprobin{t_0,x_0}{0 < x_s < \tx(s) \;\forall s\in[t_0,t], 
\tau^0=\infty} 
&\leqs \Bigprobin{t_0,x_0}{0 < x^0_s < \tx(s) \;\forall s\in[t_0,t]} \\
&\leqs \frac{\tx(t)\sqrt{a_0(t)}}{\sqrt\pi\sigma} 
\frac{\e^{-\kappa\alpha(t,t_0)/\eps}}
{\sqrt{1-\e^{-2\kappa\alpha(t,t_0)/\eps}}}.
\end{split}
\end{equation}
\end{lemma}
\goodbreak
\begin{proof}\hfill
\begin{enum}
\item   Let $g(x,s) = f(x,s) - a_0(s)x$. By assumption, $g(x,s)$ is
non-negative for $(x,s)\in\cD^+$. 
The difference $z_s = x_s - x^0_s$ satisfies the equation
\begin{equation}
\label{p:ez1}
z_s = z_{t_0} + \frac1\eps \int_{t_0}^s \bigbrak{g(x_u,u) + a_0(u)z_u} \6u
\end{equation}
with $z_{t_0}=0$. Since $g(x_s,s)\geqs0$ for $t_0\leqs s\leqs
\tau_{{\cD}^+} \wedge t$, 
\begin{equation}
\label{p:ez3}
z_s \geqs z_{t_0} + \frac1\eps \int_0^s a_0(u) z_u \6u, 
\end{equation}
follows for all such $s$ and, therefore, Gronwall's lemma yields
\begin{equation}
\label{p:ez5}
z_s \geqs z_{t_0} \e^{\kappa\alpha(s,t_0)/\eps} = 0 
\qquad \text{for all $s\in[t_0,\tau_{{\cD}^+}\wedge t]$}.
\end{equation}
This shows $x_s\geqs x^0_s$ for those $s$. Now assume
$\tau_{{\cD}^+}=\infty$ and $\tau^0=\infty$. Then, \eqref{p:ez5} 
implies that $0 < x^0_s\leqs x_s < \tx(s)$ for all $s\leqs t$,
which shows the first inequality in \eqref{ez3}.

\item   $x^0_s$ being distributed according to a normal law, we have
\begin{equation}
\label{p:ez7}
\begin{split}
\bigprobin{t_0,x_0}{0 < x^0_s < \bar x(s) \ \forall s\in[t_0,t]} 
&\leqs \bigprobin{t_0,x_0}{0 < x^0_t < \tx(t)} \\
&\leqs \frac{\tx(t)}{\sqrt{2\pi v_0(t,t_0)}},
\end{split}
\end{equation}
where the variance $v_0(t,t_0)$  
can be estimated by Lemma \ref{l_unst_var}. This proves the second
inequality in \eqref{ez3}.\qed
\end{enum}
\renewcommand{\qed}{}
\end{proof}

The previous lemma is useful only if we can control the probability that
the solution $x^0_t$ of the linearized equation returns to zero. The
following result estimates this probability and its density.

\begin{lemma}
\label{l_rz}
Let $t_0\geqs\sqrt\eps$ and assume that $x^0_{t_0} = \rho >
\sigma/\sqrt{a_0(t_0)}$. Denote by $\tau^0$ the time of the first
return of $x^0_t$ to zero. Then we have 
\begin{align}
\label{rz2a}
\probin{t_0,\rho}{\tau^0< t} &\leqs \probin{t_0,\rho}{\tau^0<
\infty} \leqs \e^{-a_0(t_0)\rho^2/\sigma^2} \\
\label{rz2b}
\dtot{}{t} \probin{t_0,\rho}{\tau^0< t} &\leqs 
\frac{2}{\sqrt\pi} \sqrt{a_0(t_0)}\,\frac\rho\sigma
\e^{-a_0(t_0)\rho^2/\sigma^2} \frac1\eps \sqrt{a_0(t)a_0(t_0)}  
\,\frac{\e^{-2\kappa\alpha(t,t_0)/\eps}}
{\sqrt{1-\e^{-2\kappa\alpha(t,t_0)/\eps}}}.
\end{align}
\end{lemma}
\goodbreak
\begin{proof}\hfill
\begin{enum}
\item   Since by symmetry, $\probin{\tau^0,0}{x^0_t\geqs0} = \frac12$ on
$\set{\tau^0< t}$, we have by the strong Markov property  
\begin{equation}
\label{prz:1}
\pcondin{t_0,\rho}{x^0_t\geqs0}{\tau^0< t} = \frac12.
\end{equation}
We now observe that 
\begin{equation}
\label{prz:2}
\begin{split}
\probin{t_0,\rho}{x^0_t\geqs0} 
& = \probin{t_0,\rho}{x^0_t\geqs0,\tau^0\geqs t} +
\probin{t_0,\rho}{x^0_t\geqs0,\tau^0< t} \\
& = \probin{t_0,\rho}{\tau^0\geqs t} +
\pcondin{t_0,\rho}{x^0_t\geqs0}{\tau^0< t}\probin{t_0,\rho}{\tau^0< t}
\\
&= 1-\probin{t_0,\rho}{\tau^0< t} + 
\tfrac12\probin{t_0,\rho}{\tau^0< t}
\\
&= 1-\tfrac12\probin{t_0,\rho}{\tau^0< t},
\end{split}
\end{equation}
which implies
\begin{equation}
\label{prz:3}
\probin{t_0,\rho}{\tau^0< t} =
2\bigbrak{1-\probin{t_0,\rho}{x^0_t\geqs0}} = 2\probin{t_0,\rho}{x^0_t< 0}.
\end{equation}

\item
Next, we use that $x^0_t$ is a Gaussian random variable with
mean $\rho\e^{\kappa\alpha(t,t_0)/\eps}$ and variance 
\begin{equation}
\label{prz:4}
v_0(t,t_0) = \frac{\sigma^2}{\eps} 
\int_{t_0}^t \e^{2\kappa\alpha(t,s)/\eps} \6 s. 
\end{equation}
By Lemma \ref{l_unst_var},
\begin{equation}
\label{prz:7}
\Xi = \frac{\rho^2\e^{2\kappa\alpha(t,t_0)/\eps}}{2v_0(t,t_0)} \geqs
a_0(t_0)\frac{\rho^2}{\sigma^2},
\end{equation}
and we thus have 
\begin{align}
\nonumber
\probin{t_0,\rho}{x^0_t< 0} 
&= \frac1{\sqrt{2\pi v_0(t,t_0)}} \int_{-\infty}^0 
\exp\Bigset{-\frac{(x-\rho\e^{\kappa\alpha(t,t_0)/\eps})^2}{2v_0(t,t_0)}}
 \6 x \\
\label{prz:8}
&= \frac1{\sqrt{2\pi}} 
\int_{-\infty}^{-\frac{\rho\e^{\kappa\alpha(t,t_0)/\eps}}{\sqrt{v_0(t,t_0)}}} 
\e^{-y^2/2} \6 y \leqs \frac12 \e^{-\Xi},
\end{align}
which proves \eqref{rz2a}, using \eqref{prz:3} and \eqref{prz:7}. 

\item
In order to compute the derivative of $\probin{t_0,\rho}{x^0_t< 0}$, 
we first note that
\begin{equation}
\label{prz:9}
\dtot{}{t} v_0(t,t_0) = \frac{\sigma^2}{\eps} +
\frac{2a_0(t)}{\eps}v_0(t,t_0).
\end{equation}
Differentiating the second line of \eqref{prz:8}, we get
\begin{align}
\nonumber
\dtot{}{t}\probin{t_0,\rho}{x^0_t< 0}
&= \frac1{\sqrt{2\pi}}
\exp\biggset{-\frac{\rho^2\e^{2\kappa\alpha(t,t_0)/\eps}}{2v_0(t,t_0)}} 
\dtot{}{t}
\biggbrak{-\frac{\rho\e^{\kappa\alpha(t,t_0)/\eps}}{\sqrt{v_0(t,t_0)}}} \\
\nonumber
&= \frac1{\sqrt{2\pi}} \e^{-\Xi}
\frac{\rho}2 \frac{\sigma^2}{\eps}
\frac{\e^{\kappa\alpha(t,t_0)/\eps}}{v_0(t,t_0)^{3/2}} \\
\label{prz:10}
&= \frac1{\sqrt{2\pi}} \frac1\rho \frac{\sigma^2}{\eps}
\frac{\e^{-\kappa\alpha(t,t_0)/\eps}}{\sqrt{v_0(t,t_0)}} \Xi \e^{-\Xi} \\
\nonumber
&\leqs \frac1{\sqrt{2\pi}} \sqrt{a_0(t_0)}\,\frac\rho\sigma
\e^{-a_0(t_0)\rho^2/\sigma^2} \frac1\eps \sqrt{2a_0(t)a_0(t_0)} 
\,\frac{\e^{-2\kappa\alpha(t,t_0)/\eps}}
{\sqrt{1-\e^{-2\kappa\alpha(t,t_0)/\eps}}},
\end{align}
where we have used the facts that $\Xi>a_0(t_0)\rho^2/\sigma^2>1$ and
that $\Xi\e^{-\Xi}$ is decreasing for $\Xi>1$. Now, \eqref{rz2b}
follows from~\eqref{prz:3}.
\qed
\end{enum}
\renewcommand{\qed}{}
\end{proof}

Assume for the moment that $x^0_t$ starts \lq\lq on the border\rq\rq\
of $\cS$, i.e. in $\rho(t_0)=h/\sqrt{a(t_0)}=\sqrt\kappa
h/\sqrt{a_0(t_0)}$. Then, by our choice
$h=2\sigma\sqrt{\abs{\log\sigma}}$, Estimate~\eqref{rz2a} shows that
the probability for $x^0_t$ to return to zero cannot exceed 
$\e^{-a_0(t_0) \rho^2/\sigma^2} = \sigma^{4\kappa}$.

We are now ready to prove the main estimate on the first exit time
$\tau_\cD$, which is the most important of our results. Since the proof is
rather involved, we restate Theorem~\ref{t_escape} here for
convenience. 

\begin{prop}[{\rm Theorem~\ref{t_escape}}\kern0pt]
\label{tt_escape}
Let $t_0\geqs\sqrt\eps$ and $\abs{x_0}\leqs \tx(t_0)$. 
Then 
\begin{equation}
\label{tesc2}
\bigprobin{t_0,x_0}{\tau_{\cD}\geqs t}\leqs 
C_0 \,\tx(t)\sqrt{a(t)} \mskip2mu \frac{\abs{\log\sigma}}{\sigma} 
\biggpar{1+\frac{\alpha(t,t_0)}\eps}
\frac{\e^{-\kappa \alpha(t,t_0)/\eps}}{\sqrt{1-\e^{-2\kappa\alpha(t,t_0)/\eps}}},
\end{equation}
where $C_0>0$ is a (numerical) constant.
\end{prop}

The strategy of the proof can be summarized as follows. The paths are
likely to leave $\cS$ after a short time. Then there are two
possibilities. Either the solution $x^0_t$ of the linear equation~\eqref{pe4}
does not return to zero, and Lemma~\ref{l_ez} shows that $x_t$ is
likely to leave $\cD$ as well. Or $x^0_t$ does return to zero. Using
the (strong) Markov property and integrating over the distribution of
the time of such a (first) return to zero, we obtain an integral
equation for an upper bound on the probability of remaining in
$\cD$. Finally, this integral equation is solved by iterations.

\goodbreak
\begin{proof}[{\sc Proof of Proposition~\ref{tt_escape}}]\hfill
\begin{enum}
\item   We first introduce some notations. Let 
\begin{equation}
\label{p:tesc1}
\Phi_t(s,x) = \bigprobin{s,x}{\tau_{\cD}\geqs t} 
= \Bigprobin{s,x}{\sup_{s\leqs u\leqs t} \frac{\abs{x_u}}{\tx(u)}< 1},
\end{equation}
and define $\rho(t) = h/\sqrt{a(t)}$. We may assume that $\rho(t)\leqs
\tx(t)$ for all $t$ (otherwise we replace $\tx$ by its maximum with
$\rho$). For $t\geqs s\geqs \sqrt\eps$ we define the quantities
\begin{align}
\label{p:tesc2}
q_t(s) &= \sup_{\abs{x}\leqs \rho(s)} \Phi_t(s,x), \\
\label{p:tesc3}
Q_t(s) &= \sup_{\rho(s)\leqs \abs{x}\leqs \tx(s)} \Phi_t(s,x).
\end{align}

\item   Let us first consider the case $\abs{x}\leqs \rho(s)$.  Recall
that  $\cS=\setsuch{(x,t)}{\abs{x} < \rho(t)}$. By Proposition~\ref{l_escape}
and the strong Markov property, we have the estimate
\begin{align}
\nonumber
\Phi_t(s,x) &= 
\bigprobin{s,x}{\tau_{\cS}\geqs t} 
+ \Bigprobin{s,x}{\tau_{\cS}<t,\sup_{\tau_{\cS}\leqs u\leqs t}
\frac{\abs{x_u}}{\tx(u)}< 1}\\
\nonumber
&\leqs \Bigpar{\frac{h}{\sigma}}^2 \e^{-\kappa\alpha(t,s)/\eps} 
+ \Bigexpecin{s,x}{\indexfct{\tau_{\cS}<t}
\Bigprobin{\tau_{\cS},x_{\tau_{\cS}}}{\sup_{\tau_{\cS}\leqs u\leqs t} 
\frac{\abs{x_u}}{\tx(u)}< 1}} \\
\label{p:tesc4}
&\leqs \Bigpar{\frac{h}{\sigma}}^2 \e^{-\kappa\alpha(t,s)/\eps} 
+ \bigexpecin{s,x}{\indicator{[s,t)}(\tau_{\cS}) Q_t(\tau_{\cS})}.
\end{align}
The second term can be estimated by integration by parts, see
Lemma~\ref{l_partint}.  Let $\bQ_t(u)$ be any upper bound on $Q_t(u)$
satisfying the hypotheses on $g$ in that lemma. Since $Q_t(u)\leqs Q_t(t) =
1$, we may assume that $\bQ_t(t)= 1$. Application of \eqref{ip1} with
$G(u)=1-(h/\sigma)^2\e^{-\kappa\alpha(u,s)/\eps}$  shows that the second
term in \eqref{p:tesc4} is bounded by
\begin{equation}
\label{p:tesc5}
\Bigpar{\frac{h}{\sigma}}^2 \e^{-\kappa\alpha(t,s)/\eps}
+ \kappa\Bigpar{\frac{h}{\sigma}}^2
\int_s^t \bQ_t(u) \frac{a(u)}{\eps}
\e^{-\kappa\alpha(u,s)/\eps} \6 u.
\end{equation}
We have thus obtained the inequality
\begin{equation}
\label{p:tesc6}
q_t(s) \leqs
2\Bigpar{\frac{h}{\sigma}}^2 \e^{-\kappa\alpha(t,s)/\eps} 
+ \kappa\Bigpar{\frac{h}{\sigma}}^2 \int_s^t \bQ_t(u)
\frac{a(u)}{\eps} \e^{-\kappa\alpha(u,s)/\eps} \6 u.
\end{equation}

\item   Consider now the case $\abs{x}\in[\rho(s),\tx(s)]$. Since $x\mapsto
f(x,t)$ is an odd function, $\Phi_t(s,x)=\Phi_t(s,-x)$ follows. Hence we
may assume that $x>0$. We consider the linear SDE~\eqref{pe4} with initial
condition $x^0_s=x$,  and denote by $\tau^0$ the time of the first return
of $x^0_t$ to zero. Then we have 
\begin{equation}
\label{p:tesc8}
\Phi_t(s,x) = \Bigprobin{s,x}{\tau^0\geqs t, \sup_{s\leqs u\leqs t}
\frac{\abs{x_u}}{\tx(u)} < 1} 
+ \Bigprobin{s,x}{ \tau^0< t, \sup_{s\leqs u\leqs t}
\frac{\abs{x_u}}{\tx(u)} < 1},
\end{equation}
and Lemma \ref{l_ez} yields
\begin{equation}
\label{p:tesc9}
\Bigprobin{s,x}{\tau^0\geqs t,\sup_{s\leqs u\leqs t}
\frac{\abs{x_u}}{\tx(u)} < 1} \leqs \frac{\tx(t)\sqrt{\kappa
a(t)}}{\sqrt\pi \sigma}
\frac{\e^{-\kappa\alpha(t,s)/\eps}}
{\sqrt{1-\e^{-2\kappa\alpha(t,s)/\eps}}}. 
\end{equation}
The second term in \eqref{p:tesc8} can be estimated using the density of
the random variable $\tau^0$, for which Lemma~\ref{l_rz} gives the bound
\begin{equation}
\label{p:tesc10}
\psi_{\tau^0}(u) = \dtot{}{u} \bigprobin{s,x}{\tau^0<u} 
\leqs \frac{2\kappa^{3/2}}{\sqrt\pi} \frac{h}{\sigma} \e^{-\kappa
h^2/\sigma^2} \frac{a(u)}\eps \frac{\e^{-2\kappa\alpha(u,s)/\eps}}
{\sqrt{1-\e^{-2\kappa\alpha(u,s)/\eps}}}.
\end{equation}
We obtain
\begin{align}
\nonumber
\Bigprobin{s,x}{\tau^0< t, \sup_{s\leqs u\leqs t}
\frac{\abs{x_u}}{\tx(u)} < 1} 
&\leqs \Bigexpecin{s,x}{\indexfct{\tau^0< t} 
\Bigprobin{\tau^0 ,x_{\tau^0}}{\sup_{\tau^0\leqs u\leqs t}
\frac{\abs{x_u}}{\tx(u)} < 1}} \\ 
\nonumber
&= \int_s^t \psi_{\tau^0}(u) \Phi_t(u,x_u) \6 u \\ 
\label{p:tesc11}
&\leqs \int_s^t \psi_{\tau^0}(u) \bigbrak{q_t(u)+Q_t(u)} \6 u.
\end{align}

\item
Before inserting the estimate \eqref{p:tesc6} for $q_t(u)$, we shall 
introduce some notations and provide bounds for certain integrals
needed in the sequel. Let
\begin{equation}
\label{p:tesc14}
g(t,s) = \frac{\e^{-\kappa\alpha(t,s)/\eps}}
{\sqrt{1-\e^{-2\kappa\alpha(t,s)/\eps}}} 
\end{equation}
and $\phi = \e^{-\kappa\alpha(t,s)/\eps}$. Then 
\begin{align}
\label{p:tesc17a}
&{}\int_s^t \frac{a(u)}\eps \e^{-\kappa\alpha(u,s)/\eps} g(u,s)\6 u 
\leqs \int_s^t \frac{a(u)}\eps g(u,s)\6 u 
\leqs \frac\pi{2\kappa} \leqs \frac2\kappa  \\ 
\label{p:tesc17b}
&{}\int_s^t \frac{a(u)}\eps \e^{-\kappa\alpha(u,s)/\eps} g(t,u)g(u,s)\6 u 
= \frac{\phi}{2\kappa} \int_0^1 \frac{\6 x}{\sqrt{x(1-x)}} 
= \frac\pi{2\kappa}\phi < \frac2\kappa\phi \\
&{}\int_s^t \frac{a(u)}\eps \e^{-\kappa\alpha(u,s)/\eps} g(t,u) \6 u 
\leqs \frac{\phi}\kappa\int_0^{\sqrt{1-\phi^2}} \frac1{1-x^2} \6 x
= \frac{\phi}\kappa \frac12
\log\frac{1+\sqrt{1-\phi^2}}{1-\sqrt{1-\phi^2}}
\nonumber \\
\label{p:tesc17c}
&\phantom{{}\int_s^t \frac{a(u)}\eps \e^{-\kappa\alpha(u,s)/\eps} g(t,u) \6 u} 
{}\leqs{} \frac{\phi}\kappa \log\frac2{\phi}
\leqs \Bigbrak{\frac1\kappa + \frac{\alpha(t,s)}\eps} \e^{-\kappa\alpha(t,s)/\eps},
\end{align}
where we used the changes of variables $\e^{-2\kappa\alpha(u,s)/\eps} =
x(1-\phi^2)+\phi^2$ in~\eqref{p:tesc17b} and
$x^2=1-\e^{-2\kappa\alpha(t,u)/\eps}$ in~\eqref{p:tesc17c}. 

\item
Now we are ready to return to our estimate on 
$\int_s^t \psi_{\tau^0}(u) q_t(u) \6 u$,
compare~\eqref{p:tesc11}. Inserting the bound~\eqref{p:tesc6} on
$q_t(u)$ yields two summands, the first one being
\begin{align}
\nonumber 
&2\Bigpar{\frac{h}{\sigma}}^2 \int_s^t \psi_{\tau^0}(u)
\e^{-\kappa\alpha(t,u)/\eps} \6 u   \\
\nonumber 
&\qquad\qquad{}\leqs \frac{4\kappa^{3/2}}{\sqrt\pi} \Bigpar{\frac
h\sigma}^{3} \e^{-\kappa h^2/\sigma^2} \int_s^t \frac{a(u)}\eps
\frac{\e^{-2\kappa\alpha(u,s)/\eps}}
{\sqrt{1-\e^{-2\kappa\alpha(u,s)/\eps}}} \e^{-\kappa\alpha(t,u)/\eps}
\6 u \\ 
\label{p:tesc12}
&\qquad\qquad{}\leqs 2\sqrt{\pi\kappa} \Bigpar{\frac h\sigma}^{3} 
\e^{-\kappa h^2/\sigma^2} \e^{-\kappa\alpha(t,s)/\eps},
\end{align}
where we used~\eqref{p:tesc17a} to bound the integral. The second summand is 
\begin{align}
\nonumber
&\kappa\Bigpar{\frac{h}{\sigma}}^2 \int_s^t \psi_{\tau^0}(u)
\int_u^t \bQ_t(v) \frac{a(v)}{\eps}
\e^{-\kappa\alpha(v,u)/\eps} \6 v\6 u \\ 
\label{p:tesc13}
&\qquad\qquad{}\leqs \kappa\sqrt{\pi\kappa}\Bigpar{\frac h\sigma}^{3}
\e^{-\kappa h^2/\sigma^2}   
\int_s^t \bQ_t(v) \frac{a(v)}{\eps} \e^{-\kappa\alpha(v,s)/\eps} \6 v,
\end{align}
where we used~\eqref{p:tesc17a} again. 

We can now collect terms. Introducing the abbreviations
\begin{equation}
\label{p:tesc16}
C = \max\Bigset{\frac{\tx(t)\sqrt{\kappa a(t)}}{\sqrt\pi \sigma},1} 
\qquad\text{and}\qquad
c = \sqrt{\pi\kappa} \Bigpar{\frac h\sigma}^{3} 
\e^{-\kappa h^2/\sigma^2}, 
\end{equation}
the previous inequalities imply that
\begin{equation}
\label{p:tesc15}
Q_t(s) \leqs C g(t,s) + c\e^{-\kappa\alpha(t,s)/\eps} 
+ c\int_s^t \bQ_t(u) \frac{a(u)}\eps \e^{-\kappa\alpha(u,s)/\eps}
\bigbrak{1+g(u,s)} \6 u.
\end{equation}

\item   
We will now iterate the bounds on $Q_t(s)$. This will show the
existence of two series $\set{a_n}_{n\geqs1}$ and
$\set{b_n}_{n\geqs1}$ such that
\begin{equation}
\label{p:tesc18}
Q_t(s) \leqs C g(t,s) + a_n \e^{-\kappa\alpha(t,s)/\eps} + b_n
\qquad \forall n.
\end{equation}
To do so, we need to assume that
\begin{equation}
\label{I_hate_it}
c \Bigpar{\frac{\alpha(T,t_0)}\eps +\frac2\kappa} 
= \sqrt{\pi\kappa}\Bigpar{\frac{\alpha(T,t_0)}\eps +\frac2\kappa}
\Bigpar{\frac h\sigma}^{3} \e^{-\kappa h^2/\sigma^2}\leqs \frac12.
\end{equation}
By our choice~\eqref{pe2b} of $h$, this condition reduces to
\begin{equation}
\label{Its_not_so_bad}
\sigma^{2\kappa}\abs{\log\sigma}^{3/4} = \Order{\sqrt\eps},
\end{equation}
which is satisfied for small enough $\eps$ by our assumption~\eqref{pe3}
on $\sigma$, provided $\kappa>1/2$. 

Using the trivial bound $\bQ_t(u)=1$ in \eqref{p:tesc15}, we find that
\eqref{p:tesc18} holds with $a_1=c$ and $b_1=3c/\kappa$. Inserting
\eqref{p:tesc18} into \eqref{p:tesc15} again, we get 
\begin{align}
\nonumber
Q_t(s) \leqs{} & C g(t,s) + c\e^{-\kappa\alpha(t,s)/\eps} \\ 
\nonumber
& + c\int_s^t \Bigbrak{C g(t,u) + a_n \e^{-\kappa\alpha(t,u)/\eps}+ b_n}
\frac{a(u)}\eps \e^{-\kappa\alpha(u,s)/\eps} \bigbrak{1+g(u,s)} \6 u \\
\leqs{} & C g(t,s)
+ c\biggbrak{1+C\Bigpar{\frac{\alpha(t,s)}\eps + \frac3\kappa} 
+ a_n\Bigpar{\frac{\alpha(t,s)}\eps + \frac2\kappa}} \e^{-\kappa\alpha(t,s)/\eps} 
+ \frac {3c}{\kappa} b_n.
\nonumber
\end{align}
By induction, we find 
\begin{align}
\nonumber
a_{n+1} 
&= c\Bigbrak{1+C\Bigpar{\frac{\alpha(t,s)}\eps + \frac3\kappa}} 
\sum_{j=0}^{n-1} \Bigbrak{c\Bigpar{\frac{\alpha(t,s)}\eps + \frac2\kappa}}^j 
+ c\Bigbrak{c\Bigpar{\frac{\alpha(t,s)}\eps + \frac2\kappa}}^n \\
\label{p:tesc20a}
&\leqs \Bigbrak{1+C\Bigpar{\frac{\alpha(t,s)}\eps + \frac3\kappa}} 
\frac c{1-c\bigpar{\frac{\alpha(t,s)}\eps + \frac2\kappa}} \\
\label{p:tesc20b}
b_{n+1} &= \Bigpar{\frac {3c}{\kappa}}^{n+1}
\end{align}
as a possible choice, where we have used the fact that $c(\alpha(t,s)/\eps
+ 2/\kappa) \leqs \frac12$ by the hypothesis \eqref{I_hate_it}. Taking the
limit $n\to\infty$, and using $c\leqs\frac\kappa4\leqs\frac14$, we
obtain  
\begin{equation}
\label{p:tesc21}
Q_t(s) \leqs C g(t,s) + \frac12 \bigpar{1+3C} \e^{-\kappa\alpha(t,s)/\eps} 
\leqs 3C g(t,s).
\end{equation}
In order to obtain also a bound on $q_t(s)$, we insert the above bound
on $Q_t(s)$ into \eqref{p:tesc6}, which yields
\begin{align}
\nonumber
q_t(s) & {}\leqs 2\Bigpar{\frac{h}{\sigma}}^2 \e^{-\kappa\alpha(t,s)/\eps}
+ 3\kappa C \Bigpar{\frac{h}{\sigma}}^2 \int_s^t \frac{a(u)}\eps
\e^{-\kappa\alpha(u,s)/\eps} g(t,u) \6 u\\
\label{p:tesc24}
& {}\leqs \Bigbrak{2+3 \kappa C \Bigpar{\frac1\kappa+
\frac{\alpha(t,s)}\eps}} \Bigpar{\frac{h}{\sigma}}^2
\e^{-\kappa\alpha(t,s)/\eps} 
\end{align}
by~\eqref{p:tesc17c}. This proves the proposition, and therefore
Theorem~\ref{\t_escape}, by taking the sum of the above estimates on
$q_t(s)$ and $Q_t(s)$. 
\qed
\end{enum}
\renewcommand{\qed}{}
\end{proof}


\subsection{Approach to $x^\star(t)$}
\label{ssec_papp}

We finally turn to the behaviour after the time $\tau=\tau_\cD > \sqrt\eps$,
when $x_t$ leaves the set $\cD$. By symmetry, we can restrict the analysis
to the case $x_\tau = \tx(\tau)$. Our aim is to prove that with high
probability, $x_t$ soon reaches a neighbourhood of $x^\star(t)$.

We start by analysing the solution $\xdettau{t}$ of the deterministic
equation 
\begin{equation}
\label{pa1}
\eps \dtot xt = f(x,t)
\end{equation}
with initial condition $\xdettau\tau = \tx(\tau)$.

\begin{prop}
\label{p_app1}
For sufficiently small $\eps$ and $T$, 
\begin{align}
\label{pa2a}
&\tx(t) \leqs \xdettau t \leqs x^\star(t) \\
\label{pa2b}
&0 \leqs x^\star(t)-\xdettau t \leqs 
C \biggbrak{\frac\eps{t^{3/2}} + \bigpar{x^\star(\tau)-\tx(\tau)}
\e^{-\mm\alpha(t,\tau)/\eps}} \\
\label{pa2c}
&0 \leqs \xdetin{\sqrt\eps}t - \xdettau t \leqs 
\bigpar{\xdetin{\sqrt\eps}\tau - \tx(\tau)} \e^{-\mm\alpha(t,\tau)/\eps}
\end{align}
for all $t\in[\tau,T]$ and all $\tau\in[\sqrt\eps,T]$, where $C>0$ is a
constant depending only on $f$.
\end{prop}
\goodbreak
\begin{proof}\hfill
\begin{enum}
\item   Whenever $\xdettau{t}=x^\star(t)$, we have 
\begin{equation}
\label{p:app1:1}
\eps\dtot{}{t} \bigpar{x^\star(t)-\xdettau{t}} = \eps\dtot{x^\star(t)}t - 
f(x^\star(t),t) = \eps\dtot{x^\star(t)}t \geqs 0,
\end{equation}
which shows that $\xdettau{t}$ can never become larger than $x^\star(t)$. 
Similarly, whenever $\xdettau{t}=\tx(t)$, we get
\begin{equation}
\label{p:app1:2}
\begin{split}
\eps\dtot{}{t} \bigpar{\xdettau{t}-\tx(t)} &= 
f(\tx(t),t) - \eps\dtot{\tx(t)}t \\
&= \sqrt\lambda\,(1-\lambda) t^{3/2} \bigbrak{1+\orderone{T}} 
- \eps\frac{\sqrt\lambda}{2\sqrt t} \bigbrak{1+\orderone{T}} > 0
\end{split}
\end{equation}
provided $\lambda<\frac12[1 - \orderone{T}]$, which shows that
$\xdettau{t}$ can never become smaller than $\tx(t)$. This completes
the proof of~\eqref{pa2a}. 

\item   We now introduce the difference $\ydettau t = x^\star(t) - \xdettau
t$. Using Taylor's formula, one immediately obtains that $\ydettau t$
satisfies the ODE
\begin{equation}
\label{p:app1:3}
\eps\dtot yt = a^\star(t)y + b^\star(y,t) + \eps x^{\star\prime}(t)
\end{equation}
where
\begin{equation}
\label{p:app1:4}
\begin{split}
a^\star(t) &\leqs -a_0^\star t \\
0 \leqs b^\star(y,t) &\leqs M^\star\sqrt t \,y^2 \\
\xstarprime(t) &\leqs \frac{K^\star}{\sqrt t},
\end{split}
\end{equation}
with $a_0^\star = 2[1+\orderone{T}]$, $M^\star = 3[1+\orderone{T}]$ and
$K^\star = \frac12[1+\orderone{T}]$. We first consider the particular
solution $\yhatdet_t$ of \eqref{p:app1:3} starting at time
$4\sqrt\eps$ in \smash{$\yhatdet_{4\sqrt\eps}=0$}. 
By~\eqref{p:app1:1}, we know that $\yhatdet_t\geqs 0$ for all
$t\geqs4\sqrt\eps$. We will use the fact that 
\begin{equation}
\label{p:app1:5}
\begin{split}
\int_\tau^t  \frac1{\sqrt s}\e^{-a_0^\star(t^2-s^2)/4\eps} \6s 
&\leqs \int_\tau^t \frac1{\sqrt s}\e^{-a_0^\star t(t-s)/4\eps} \6s \\
&\leqs \frac{4\eps}{a_0^\star t^{3/2}}\int_0^\x\frac{\e^{-u}}{\sqrt{1-u/\x}}\6u
< c_0 \frac\eps{t^{3/2}},
\end{split}
\end{equation}
where $c_0 = 8/a_0^\star$. We have used the transformation $s=t-4\eps
u/(a_0^\star t)$, introduced $\x=a_0^\star t^2/4\eps$ and bounded the last
integral by $2$. We now introduce the first exit time 
$\hat\tau = \inf\setsuch{t\geqs4\sqrt\eps}{\yhatdet_t\geqs c_0\eps
t^{-3/2}}$. For $4\sqrt\eps \leqs t \leqs \hat\tau$, we have
\begin{equation}
\label{p:app1:6}
a^\star(t)y + b^\star(y,t) 
\leqs \Bigpar{-a_0^\star\,t + M^\star\sqrt t\,c_0\frac{\eps}{t^{3/2}}} y 
\leqs -a_0^\star 
\Bigpar{1 - \frac{c_0M^\star}{16a_0^\star}} t y. 
\end{equation}
Since $M^\star/(a_0^\star)^2 = \frac34 [1+ \orderone{}]$, the term in
brackets can be assumed to be larger than $\frac12$. Hence
\eqref{p:app1:3} shows that 
\begin{equation}
\label{p:app1:7}
\eps\dtot{\yhatdet}{t} 
\leqs -\frac{a_0^\star}2 \mskip.5mu t \mskip1mu \yhatdet 
+ \eps\frac{K^\star}{\sqrt t},
\end{equation}
which implies
\begin{equation}
\label{p:app1:7b}
\yhatdet_t \leqs K^\star 
\int_\tau^t  \frac{\e^{-a_0^\star(t^2-s^2)/4\eps}}{\sqrt s} \6s 
< K^\star c_0 \frac{\eps}{t^{3/2}}.
\end{equation}
Since $K^\star = \frac12[1 +\orderone{}]$, we obtain $\yhatdet_t<c_0\eps
t^{-3/2}$, and thus $\hat\tau=\infty$. This shows
\begin{equation}
\label{p:app1:8}
0 \leqs \yhatdet_t \leqs K^\star c_0 \frac{\eps}{t^{3/2}} 
\qquad \text{for $4\sqrt\eps \leqs t \leqs T$.}
\end{equation}

\item   Let $\tau\geqs\sqrt\eps$ and $0\leqs y_1< y_2 \leqs
x^\star(\tau)-\tx(\tau)$ be given. Let $y^{(1)}_t$ and $y^{(2)}_t$ be
solutions of \eqref{p:app1:3} with initial conditions $y^{(1)}_\tau=y_1$ and
$y^{(2)}_\tau=y_2$, respectively. Then there exists a $\theta\in[0,1]$
such that the difference $z_t=y^{(2)}_t-y^{(1)}_t$ satisfies
\begin{equation}
\label{p:app1:9}
\eps\dtot zt = -\sdpar fx (x^\star(t)-y^{(1)}_t-\theta z,t) 
\leqs -\mm a(t) z,
\end{equation}
where we have used \eqref{pa2a} and \eqref{pp16b}. It follows that 
\begin{equation}
\label{p:app1:10}
0\leqs y^{(2)}_t-y^{(1)}_t \leqs (y_2-y_1) \e^{-\mm\alpha(t,\tau)/\eps},
\end{equation}
which proves \eqref{pa2c} in particular. 
If $\tau \geqs 4\sqrt\eps$, we can use the relation 
$x^\star(t) - \xdettau t = \yhatdet_t + (\ydettau t-\yhatdet_t)$ to show
that 
\begin{equation}
\label{p:app1:11}
x^\star(t) - \xdettau t \leqs K^\star c_0 \frac\eps{t^{3/2}} + 
\bigpar{x^\star(\tau)-\tx(\tau)} \e^{-\mm \alpha(t,\tau)/\eps},  
\end{equation}
which proves \eqref{pa2b} for $\tau \geqs 4\sqrt\eps$.  Finally, if
$\sqrt\eps \leqs \tau \leqs 4\sqrt\eps$, we can use the fact that 
$x^\star(t) - \xdettau t \leqs x^\star(t) - \xdetin {4\sqrt\eps}t$ to
prove that \eqref{pa2b} holds for some constant $C>0$.
\qed
\end{enum}
\renewcommand{\qed}{}
\end{proof}

Let us now consider the process $y_t = y^\tau_t = x_t - \xdettau t$, starting
at time $\tau$ in $y_\tau=0$, which describes the deviation due to noise
from the deterministic solution $\xdettau t$. It satisfies the SDE
\begin{equation}
\label{pa3}
\6y_t = \frac1\eps \bigbrak{\atau(t)y + \btau(y_t,t)} \6t +
\frac\sigma{\sqrt\eps} \6W_t,
\end{equation}
where we have introduced 
\begin{equation}
\label{pa4}
\begin{split}
\atau(t) &= \sdpar fx(\xdettau t,t) \\
\btau(y,t) &= f(\xdettau t+y,t) - f(\xdettau t) - \atau(t)y.
\end{split}
\end{equation}
The following bounds are direct consequences of Taylor's formula and
Proposition \ref{p_app1}:
\begin{align}
\label{pa5a}
&a^\star(t) \leqs \atau(t) \leqs \ta(t) \\
\label{pa5b}
&\atau(t) = a^\star(t) + \BigOrder{\frac\eps t} +
\Order{t\e^{-\mm\alpha(t,\tau)/\eps}} \\
\label{pa5c}
&(\atau)'(t) = \BigOrder{1 + \frac{t^2}\eps \e^{-\mm\alpha(t,\tau)/\eps}} \\
\label{pa5d}
&\abs{\btau(y,t)} \leqs 3My^2 \bigpar{x^\star(t)+\abs{y}}, \qquad 
\text{valid for $x^\star(t)+\abs{y}\leqs d.$ }
\end{align}

For comparison, we will also consider the linear SDE
\begin{equation}
\label{pa6}
\6y^0_t = \frac1\eps \atau(t) y^0_t \6t + \frac\sigma{\sqrt\eps} \6W_t.
\end{equation}
Let $\alphatau(t,s) = \int_s^t \atau(u)\6u$ and denote by
\begin{equation}
\label{pa7}
v^\tau(t) = \frac{\sigma^2}\eps \int_\tau^t \e^{2\alphatau(t,s)/\eps}\6s
\end{equation}
the variance of $y^0_t$. Again we introduce and investigate a function 
\begin{equation}
\label{pa8}
\ztau(t) = \frac1{2\abs{\ta(\tau)}} \e^{2\alphatau(t,\tau)/\eps} 
+ \frac1\eps \int_\tau^t \e^{2\alphatau(t,s)/\eps}\6s.
\end{equation}

\begin{lemma}
\label{p_app2}
The function $\ztau(t)$ satisfies the following relations for $\tau\leqs
t\leqs T$:
\begin{align}
\label{pa9a}
&\ztau(t) = \frac1{2\abs{\ta(t)}} + \BigOrder{\frac\eps{t^3}} +
\BigOrder{\frac1t \e^{-\mm\alpha(t,\tau)/\eps}} \\
\label{pa9b}
&\frac1{2\abs{a^\star(t)}} \leqs \ztau(t) \leqs \frac1{2\abs{\ta(\tau)}} \\
\label{pa9c}
&(\ztau)'(t) \leqs \frac1\eps.
\end{align}
\end{lemma}
\goodbreak
\begin{proof}\hfill
\begin{enum}
\item   By integration by parts, we find
\begin{equation}
\label{p:app2:1}
\ztau(t) = \frac1{2\abs{\ta(t)}} - \frac12 \int_\tau^t
\frac{(\atau)'(s)}{\atau(s)^2} \e^{2\alphatau(t,s)/\eps} \6s.
\end{equation}
The relation $\abs{\atau(s)} \geqs \abs{\ta(s)} \geqs \mm\abs{a(s)}$
together with \eqref{pa5c} yields
\begin{equation}
\label{p:app2:2}
\biggabs{\int_\tau^t
\frac{(\atau)'(s)}{\atau(s)^2} \e^{2\alphatau(t,s)/\eps} \6s} 
\leqs \text{\it const}{} \int_\tau^t \Bigpar{\frac1{s^2} + \frac1\eps
\e^{-\mm\alpha(s,\tau)/\eps}} \e^{-2\mm\alpha(t,s)/\eps} \6s.
\end{equation}
The second term in brackets gives a contribution of order
$\frac1t\e^{-\mm\alpha(t,\tau)/\eps}$. In order to estimate the contribution
of the first term, we perform the change of variables $u=\mm(t^2-s^2)/2\eps$,
thereby obtaining
\begin{equation}
\label{p:app2:3}
\int_\tau^t \frac1{s^2}\e^{-\mm(t^2-s^2)/2\eps}\6s 
= \frac\eps{\mm t^3} \int_0^{\x-\x_0} \frac{\e^{-u}}{(1-u/\x)^{3/2}} \6u
\leqs \frac\eps{\mm t^3}
\Bigbrak{2^{3/2} + 2 \frac{\x^{3/2}\e^{-\x/2}}{\sqrt{\x_0}}},
\end{equation}
where $\x=\mm t^2/2\eps$ and $\x_0 = \mm\tau^2/2\eps$. The last inequality
is obtained by splitting the integral at $\x/2$. Using the fact that
$t^3\e^{-\mm t^2/4\eps} \leqs (6\eps/\mm)^{3/2}\e^{-3/2}$ for all $t\geqs
0$, we reach the conclusion that this integral is bounded by a
constant times $\eps/t^3$, which completes the proof of \eqref{pa9a}. 

\item   We now use the fact that $\ztau(t)$ solves the ODE
\begin{equation}
\label{p:app2:4}
\dtot{\ztau}{t} = \frac1\eps \bigpar{2\atau(t)\ztau + 1}, 
\qquad
\ztau(\tau) = \frac1{2\abs{\ta(\tau)}}.
\end{equation}
Then, \eqref{pa9c} is an immediate consequence of this relation,
and~\eqref{pa9b} is obtained from the fact that 
\begin{equation}
\label{p:app2:5}
\dtot{\ztau(t)}t 
= \frac1\eps \Bigpar{-\frac{\abs{\atau(t)}}{\abs{\ta(\tau)}}+1} \leqs 0,
\end{equation}
whenever $\ztau(t)=1/2\abs{\ta(\tau)}$, and 
\begin{equation}
\label{p:app2:6}
\dtot{}t \Bigpar{\ztau(t) - \frac1{2\abs{a^\star(t)}}}
= \frac1\eps \Bigpar{-\frac{\abs{\atau(t)}}{\abs{a^\star(t)}}+1} - 
\frac{{a^\star}'(t)}{2a^\star(t)^2} \geqs 0,  
\end{equation}
whenever $\ztau(t)=1/2\abs{a^\star(\tau)}$. Here we used~\eqref{pa5a}
and the monotonicity of $\ta(t)$ for small $t$. 
\qed
\end{enum}
\renewcommand{\qed}{}
\end{proof}

We note that Lemma~\ref{p_app2} and the bounds \eqref{pa5a} on $\atau$
imply the existence of constants $c_+\geqs c_->0$, depending only on $f$ and
$T$, such that 
\begin{equation}
\label{pa10}
\frac{c_-}t \leqs \ztau(t) \leqs \frac{c_+}t 
\qquad \forall t\in[\tau,T].
\end{equation}
We can now easily prove that $y^0_t$ remains in a strip of width
$h\sqrt{\ztau}$ with high probability, in much the same way as in
Proposition~\ref{l_ns3}.

\begin{prop}
\label{l_pa1}
For sufficiently small $T$ and $\eps$, and all $t\in[\tau,T]$, 
\begin{equation}
\label{pa11a}
\Bigprobin{\tau,0}{\sup_{\tau\leqs s\leqs t}
\frac{\abs{y^0_s}}{\sqrt{\ztau(s)}} \geqs h} 
\leqs C^\tau(t,\eps) \exp\Bigset{-\frac12\frac{h^2}{\sigma^2} 
\bigbrak{1-r(\eps)}},
\end{equation}
where $r(\eps)=\Order{\eps}$ and 
\begin{equation}
\label{pa11b}
C^\tau(t,\eps) = \frac{\abs{\alphatau(t,\tau)}}{\eps^2} + 2.
\end{equation}
\end{prop}
\begin{proof}
Let $K=\intpartplus{\abs{\alphatau(t,\tau)}/2\eps^2}$ and define a partition
$\tau=u_0 < \dots < u_K=t$ of $[\tau,t]$ by 
\begin{equation}
\label{p:lpa1:1}
\abs{\alphatau(u_k,\tau)} = 2\eps^2 k, 
\qquad k=1,\dots,K-1.
\end{equation}
Since $\atau(s)\leqs \ta(s)\leqs -\mm s/2$, we obtain
$u_k-u_{k-1}\leqs 4\eps^2/(\mm u_{k-1})$ for all $k$. Now we can
proceed as in the proof of Proposition~\ref{l_ns3}.
\end{proof}

We can now compare the solutions of the linear and the nonlinear
equation. To do so, we define the events 
\begin{align}
\label{pa12a}
\Omega_t(h) &= 
\bigsetsuch{\w}{\abs{y^\tau_s} < h\sqrt{\ztau(s)} \; \forall s\in[\tau,t]}
\\
\label{pa12b}
\Omega^0_t(h) &= 
\bigsetsuch{\w}{\abs{y^0_s} < h\sqrt{\ztau(s)} \; \forall s\in[\tau,t]}.
\end{align}
The following proposition shows that $y^\tau_t$ and $y^0_t$ differ only
slightly. 

\begin{prop}
\label{l_pa2}
Let $\gamma = 1 \,\vee\, 48M(2+\sqrt{c_+})c_+^2/\sqrt{c_-}$ and assume
$h<\tau/\gamma$ as well as $h\leqs
[d-x^\star(t)]\sqrt\tau/(2\sqrt{c_+})$. Then 
\begin{align}
\label{pa13a}
\Omega_t(h) &\subas \Omega^0_t \Bigpar{\Bigbrak{1+\gamma\frac
h\tau} h} \\
\label{pa13b}
\Omega^0_t(h) &\subas \Omega_t \Bigpar{\Bigbrak{1+\gamma\frac
h\tau} h}.
\end{align}
\end{prop}
\goodbreak
\begin{proof}
Assume first that $\w\in\Omega^0_t(h)$. We introduce the difference
$z_s=y^\tau_s-y^0_s$, set $\delta=\gamma h/\tau < 1$, and define the
first exit time 
\begin{equation}
\label{p:lpa2:2}
\hat\tau = \inf\bigsetsuch{s\in[\tau,t]} 
{\abs{z_s}\geqs\delta h \sqrt{\ztau(s)}\,} 
\in [\tau,t]\cup\set{\infty}.
\end{equation}
On $A = \Omega^0_t(h) \cap \set{\hat\tau<\infty}$, we get by the
estimate~\eqref{pa5d} on $\btau$, Lemma~\ref{p_app2} and \eqref{pa10}
\begin{align}
\nonumber
\abs{z_s} &\leqs \frac1\eps \int_\tau^t \e^{\alphatau(s,u)/\eps}
\abs{\btau(y_u,u)} \6u \\
\label{p:lpa2:1}
&\leqs 6M(1+\delta)^2 \Bigpar{2c_+\frac h\tau + (1+\delta)c_+^{3/2} 
\frac{h^2}{\tau^2}}
\frac{c_+}{\sqrt{c_-}} \,h\sqrt{\ztau(s)} < \delta h \sqrt{\ztau(s)},
\end{align}
for all $s\in[\tau,\hat\tau]$, which leads to a contradiction for
$s=\hat\tau$. We conclude that $\fP(A)=0$ and thus $\abs{z_s}\leqs
\gamma h^2\sqrt{\ztau(s)}/\tau$ for all $s$ in $[\tau,t]$, which proves
\eqref{pa13b}. The inclusion \eqref{pa13a} is a straightforward consequence
of the same estimates.
\end{proof}

Now, the following corollary is  a direct consequence of the two
preceding propositions.

\begin{cor}
\label{c_pa}
There exists $h_0$ such that if $h<h_0 \tau$, then
\begin{equation}
\label{pa14}
\biggprobin{\tau,\tx(\tau)}{\sup_{\tau\leqs s\leqs t}
\frac{\abs{x_s-\xdettau{s}}}{\sqrt{\ztau(s)}} > h} 
\leqs C^\tau(t,\eps) \exp\Bigset{-\frac12 \frac{h^2}{\sigma^2}
\Bigbrak{1-\Order{\eps}-\BigOrder{\frac h\tau}}},
\end{equation}
where $C^\tau(t,\eps)$ is given by \eqref{pa11b}. 
\end{cor}


\appendix
\section*{Appendix}
\label{appendix}

\makeatletter
\renewcommand\theequation{A.\@arabic\c@equation}
\renewcommand\thetheorem{A.\@arabic\c@theorem}
\makeatother

\setcounter{equation}{0}
\setcounter{theorem}{0}


The appendix provides two lemmas needed in Sections~\ref{sec_nobif}
and~\ref{sec_pitchfork}. The first one uses exponential martingales to
deduce an exponential bound on the probability that a stochastic integral
exceeds a given value.

\begin{lemma}
\label{p_ns1}
Let $\varphi(u)$ be a Borel-measurable deterministic function such that 
\begin{equation}
\label{ns10}
\Phi(t) = \int_0^t \varphi(u)^2 \6u
\end{equation}
exists. Then
\begin{equation}
\label{ns10a}
\Bigprob{\sup_{0\leqs s \leqs t} \int_0^s \varphi(u) \6W_u \geqs \delta} 
\leqs \exp\biggset{-\frac{\delta^2}{2\Phi(t)}}
\end{equation}
\end{lemma}
\begin{proof}
Let $P$ denote the left-hand side of \eqref{ns10a}. 
For any $\gamma>0$, we have 
\begin{equation}
\label{ns10b}
P = 
\Bigprob{\sup_{0\leqs s \leqs t} \exp\Bigset{\gamma \int_0^s
\varphi(u) \6W_u} \geqs \e^{\gamma \delta}} 
\leqs \Bigprob{\sup_{0\leqs s \leqs t} M_s
\geqs \e^{\gamma \delta - \frac{\gamma^2}{2} \Phi(t)}},
\end{equation}
where 
\begin{equation}
\label{ns10bb}
M_s = \exp\Bigset{\int_0^s \gamma \varphi(u)
\6W_u -\tfrac12 \int_0^s \gamma^2 \varphi(u)^2 \6u}
\end{equation}
is an (exponential) martingale, satisfying 
$\E\{M_t\} = \E\{M_0\} = 1$, 
which implies by Doob's submartingale inequality,
that 
\begin{equation}
\label{ns10c}
\Bigprob{\sup_{0\leqs s\leqs t}M_s\geqs\lambda} \leqs
\frac{1}{\lambda} \bigexpec{M_t} =\frac1\lambda.
\end{equation}
This gives us
\begin{equation}
\label{ns10d}
P \leqs \e^{-\gamma \delta + \frac{\gamma^2}{2} \Phi(t)},
\end{equation}
and we obtain the result by optimizing \eqref{ns10d} over $\gamma$.
\end{proof}

The following lemma allows to estimate expectation values
by integration by parts.

\begin{lemma}
\label{l_partint}
Let $\tau\geqs s_0$ be a random variable satisfying 
$F_\tau(s) = \prob{\tau<s} \geqs G(s)$ for some continuously
differentiable function $G$. Then
\begin{equation}
\label{ip1}
\bigexpec{\indicator{[s_0,t)}(\tau) g(\tau)} 
\leqs g(t) \bigbrak{F_\tau(t)-G(t)} + \int_{s_0}^t g(s)G'(s)\6 s
\end{equation}
holds for all $t>s_0$ and all functions $0\leqs g\leqs 1$
satisfying the two conditions
\begin{itemiz}
\item   there exists an $s_1\in(s_0,\infty]$ such that $g$ is
continuously differentiable and increasing on $(s_0,s_1)$; 
\item   $g(s)=1$ for all $s\geqs s_1$.
\end{itemiz}
\end{lemma}

\begin{proof}
First note that for all $t\leqs s_1$, 
\begin{align}
\nonumber
\int_{s_0}^t g'(s)\prob{\tau\geqs s}\6 s 
&= \Bigexpec{\int_{s_0}^{t\wedge\tau} g'(s)\6 s} \\
\nonumber
&= \expec{g(t\wedge\tau)} - g(s_0) \\
\label{p:ip1}
&= \expec{\indicator{[s_0,t)}(\tau) g(\tau)} + g(t)\prob{\tau\geqs t} - g(s_0)
\end{align}
which implies, by integration by parts,
\begin{align}
\nonumber
\expec{\indicator{[s_0,t)}(\tau) g(\tau)} 
&= \int_{s_0}^t g'(s) \bigbrak{1-F_\tau(s)} \6 s 
- g(t) \bigbrak{1-F_\tau(t)} + g(s_0) \\
\label{p:ip2}
&\leqs \int_{s_0}^t g(s)G'(s)\6 s + g(t) \bigbrak{F_\tau(t)-G(t)},
\end{align}
where we have used $F_\tau(s)\geqs G(s)$ and $G(s_0)\leqs F(s_0)=0$. This
proves the assertion in the case $t\leqs s_1$. In the case $t>s_1$, we have
\begin{align}
\nonumber
\expec{\indicator{[s_0,t)}(\tau) g(\tau)} 
&= \expec{\indicator{[s_0,s_1)}(\tau) g(\tau)} 
+ \prob{\tau\in[s_1,t)} \\
\nonumber
&\leqs \int_{s_0}^{s_1} g(s) G'(s)\6 s 
+ g(s_1)\bigbrak{F_\tau(s_1)-G(s_1)} 
+  \bigbrak{F_\tau(t)-F_\tau(s_1)} \\
\label{p:ip3}
&= \int_{s_0}^t g(s) G'(s)\6 s 
-  \bigbrak{G(t)-G(s_1)} 
+  \bigbrak{F_\tau(t)-G(s_1)},
\end{align}
where we have used that $g(s)=1$ holds for all $s\in[s_1,t]$. This
proves the assertion for $t>s_1$.
\end{proof}


\bigskip\bigskip\noindent
{\small 
Nils Berglund \\ 
{\sc Georgia Institute of Technology} \\ 
Atlanta, GA 30332-0430, USA \\
{\it and} \\
{\sc Weierstra\ss\ Institute for Applied Analysis and Stochastics} \\
Mohrenstra{\ss}e~39, 10117~Berlin, Germany \\
{\it E-mail address: }{\tt berglund@wias-berlin.de}

\bigskip\noindent
Barbara Gentz \\ 
{\sc Weierstra\ss\ Institute for Applied Analysis and Stochastics} \\
Mohrenstra{\ss}e~39, 10117~Berlin, Germany \\
{\it E-mail address: }{\tt gentz@wias-berlin.de}
}

\end{document}